\documentclass[openany]{book}
\usepackage[T1]{fontenc}
\usepackage[utf8]{inputenc}
\usepackage[italian]{babel}
\usepackage{amsmath}
\usepackage{amssymb}
\usepackage{xcolor}
\usepackage{MnSymbol}
\usepackage{amsthm,enumitem}
\usepackage{tikz-cd}
\usepackage{amssymb,amsmath,amsthm} 
\usepackage{bm}
\usepackage{makeidx}
\usepackage{appendix}
\usepackage{cite}

\newtheorem{thm}{Teorema}[chapter]
\newtheorem{cor}[thm]{Corollario}
\newtheorem{lem}[thm]{Lemma}
\newtheorem{prop}[thm]{Proposizione}
\theoremstyle{definition}
\newtheorem{defn}[thm]{Definizione}
\newtheorem{oss}[thm]{Osservazione}
\newtheorem{exa}[thm]{Esempio}

\newtheorem*{conj}{Congettura}

\date{}

\newcommand{\E}{{\mathbb{E}^2}}
\newcommand{\ka}{\kappa}
\newcommand{\mnk}{{M_{\ka}^n}}
\newcommand{\de}{\partial}
\newcommand{\cat}[1]{\operatorname{CAT}(#1)}
\newcommand{\catk}{\operatorname{CAT}(\ka)}
\newcommand{\mm}[1]{M_{\ka}^{#1}}
\newcommand{\sh}[1]{\operatorname{Shapes}(#1)}
\newcommand{\supp}{\operatorname{supp}}
\newcommand{\St}{\operatorname{St}}
\newcommand{\st}{\operatorname{st}}
\newcommand{\Lk}{\operatorname{Lk}}

\newcommand{\kl}{[\!\![}
\newcommand{\cay}{\operatorname{Cay}}

\makeindex
\begin{document}

\frontmatter
\thispagestyle{empty}
\clearpage
\begin{center}
	\large{\textbf{UNIVERSIT\`A DEGLI STUDI DI MILANO-BICOCCA}}\\
	\mbox{Dipartimento di Matematica e Applicazioni}\\
	\mbox{Corso di Laurea Magistrale in Matematica} 
\end{center}
\addvspace{1cm}
\begin{figure}[h]
	\centering
	\includegraphics*[width=4cm]{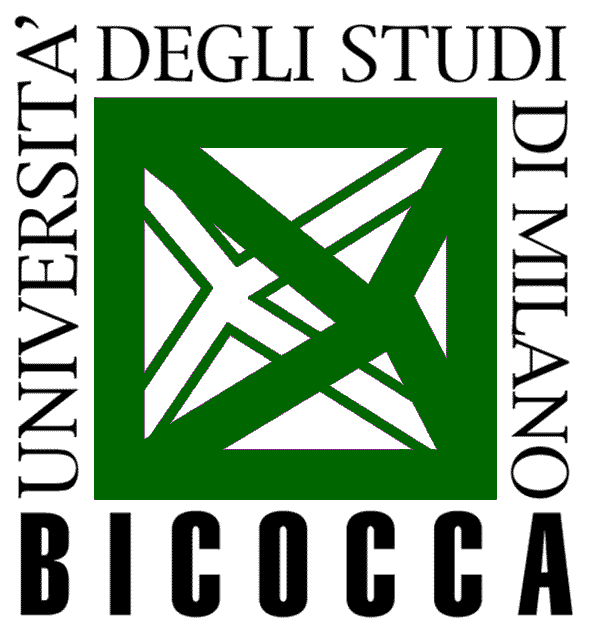}
\end{figure}

\addvspace{1cm}

\begin{center}
	\Huge
	\textbf{TEORIA GEOMETRICA DEI GRUPPI\\}
	\large
	\textbf{\Large Spazi CAT(0), Teorema di Gromov e\\}
	\textbf{\Large\textit{Oriented Right-Angled Artin Groups}\\}
	\par
\end{center}

\addvspace{1cm}
\begin{flushleft}
	\begin{tabbing}
		Relatore: Prof. Thomas S. Weigel\\
	\end{tabbing}
\end{flushleft}

\addvspace{1cm}
\begin{flushright}
	\begin{tabbing}
		\hspace{260pt}
		\= Tesi di Laurea di:\\
		\> Simone Blumer \\
		\> Matricola 805239\\
	\end{tabbing}
\end{flushright}

\begin{center}
	Anno Accademico 2019 - 2020
\end{center}
\thispagestyle{empty}
\frontmatter
	\chapter*{Introduzione}
	Lo scopo di questo elaborato finale di Laurea Magistrale è quello di introdurre la nozione di spazio a curvatura limitata dall'alto e di presentarne alcune applicazioni nel contesto della Teoria Geometrica dei Gruppi. \\
	
	La tesi consiste di due macro-sezioni, una di carattere teorico e l'altra applicativa. In particolare presentiamo un'applicazione della nozione di spazio $\cat 0$ e una del teorema di Gromov per i complessi cubici.
	Nel primo capitolo richiamiamo alcune nozioni elementari della teoria degli spazi metrici e introduciamo il concetto di spazio $\cat 0$. Uno spazio metrico geodetico è uno spazio $\cat 0$ se i triangoli in esso sono "più stretti" di quelli nel piano euclideo.\\
	Nel secondo capitolo vediamo una prima applicazione degli spazi $\cat 0$: il Teorema di Punto Fisso di Bruhat e Tits. Il teorema si colloca nello studio dei punti fissi di gruppi che agiscono su spazi. La dimostrazione è elementare, ma fornisce un buon esempio dell'utilità di spazi siffatti.\\
	Nel terzo capitolo esponiamo una generalizzazione del concetto di spazio $\cat 0$: lo spazio modello per il confronto dei triangoli non è il piano euclideo ma appartiene a una classe più ampia di spazi altamente simmetrici. A tal fine abbiamo bisogno di richiamare alcune proprietà di questi spazi modello.\\
	Nel quarto capitolo definiamo delle strutture combinatorie costruite attaccando particolari sottoinsiemi degli spazi modello, in un modo che risulta simile alla costruzione dei complessi di celle. Data la loro natura combinatorica, gli spazi costruiti forniscono una buona classe su cui studiare le azioni di alcuni gruppi finitamente generati e derivare proprietà degli stessi.\\
	Nel capitolo 5, dimostriamo il Criterio di Gromov per complessi cubici. Il risultato fornisce una equivalenza tra alcune proprietà metriche (globali) e combinatoriche (locali) dei complessi.\\
	Come applicazione del teorema sopra citato, nell'ultimo capitolo introduciamo una classe di gruppi che generalizza i ben noti gruppi di Artin (RAAG), investigandone alcune prime proprietà. Questi gruppi agiscono su un complesso cubico che ammette una metrica $\cat 0$. Costruiamo inoltre uno spazio classificante per tali gruppi e calcoliamo il loro anello di coomologia in alcuni casi particolari.  

\tableofcontents

\mainmatter

\part{}
\chapter{Spazi metrici, geodetiche e spazi CAT(0)}
La Teoria Geometrica dei Gruppi si occupa dello studio di proprietà algebriche dei gruppi tramite argomentazioni di carattere geometrico e topologico. In particolare, si  studiano le proprietà di un gruppo tramite le sue azioni su particolari spazi.  
In questa sezione richiamiamo alcune definizioni elementari necessarie per il seguito. Diamo inoltre una definizione di spazio $\cat 0$ in modo da illustrare alcune proprietà che ne derivano immediatamente. 

\section{Spazi metrici}
\begin{defn}
	Sia $X$ un insieme. Una \textit{pseudometrica} su $X$ è un'applicazione $$d:X\times X\to \mathbb R$$ tale che $\forall x,y,z\in X$: 
	\begin{enumerate}
		\item $d(x,y)=d(y,x)\geq 0$,
		\item $d(x,x)=0$,
		\item valga la \textit{disuguaglianza triangolare}: $d(x,y)\leq d(x,z)+d(z,y)$.
	\end{enumerate}
	Una pseudometrica $d$ su $X$ è una \textit{metrica} se soddisfa la seguente condizione di non-degenerazione: $\forall x,y\in X$, $d(x,y)=0$ implica $x=y$.\\
	Se $d$ è una metrica, la coppia $(X,d)$ (o anche, più semplicemente, $X$) è detta \textit{spazio metrico}. Gli spazi metrici costituiscono una buona classe di spazi su cui studiare le azioni di gruppi perché consentono di tradurre alcune proprietà dello spazio, ottenute tramite calcoli puramente geometrici, in proprietà algebriche del gruppo. 
\end{defn}

Gli spazi metrici sono \textit{spazi topologici}, i cui aperti sono i sottoinsiemi $U\subseteq X$ tali che $\forall x\in U$, esiste una palla $B_r(x)=\{y\in X\ \vert \ d(x,y)< r\}$, $r>0$, tutta contenuta in $U$. Le funzioni continue rivestono quindi un ruolo importante nello studio degli spazi metrici. Tuttavia, esiste una classe di funzioni propria degli spazi metrici. Diremo che una funzione $f:X\to X'$ tra due spazi metrici è una \textit{isometria} se $\forall x,y\in X$ vale \[d_X(x,y)=d_{X'}(f(x),f(y)).\]
Le isometrie sono chiaramente funzioni continue injettive.\\

Nel seguito, sia $X$ uno spazio metrico.
\begin{defn}
	Una \textit{geodetica} in $X$ è una funzione continua $c:[0,l]\to X$ tale che $\forall t,t'\in [0,l]$: $$d(c(t),c(t'))=\vert t-t'\vert.$$
	Il \textit{punto iniziale} (risp. \textit{terminale}) di $c$ è $c(0)$ (risp. $c(1)$); i punti iniziale e terminale si dicono \textit{estremi} di $c$, e si dice che $c$ \textit{congiunge} tali punti. Il numero reale non negativo $l=d(c(0),c(1))$ è chiaramente la distanza tra gli estremi di $c$ e si chiama \textit{lunghezza della geodetica}.\\
	L'immagine $\vert c\vert=c([0,l])$ della geodetica $c$ è detta \textit{segmento geodetico}.
\end{defn} 
Chiaramente, se una funzione continua $c:I\to X$ definita su un qualche intervallo compatto $I$ soddisfa $d(c(t),c(t'))=\lambda\vert t-t'\vert$,  per qualche $\lambda\in\mathbb R_+$, la sua immagine in $X$ è un segmento geodetico. La funzione $c$ si dice \textit{geodetica a velocità costante} $\lambda$ (oppure \textit{geodetica riparametrizzata}).
\begin{defn}
	Uno spazio $X$ è detto \textit{geodetico} se per ogni coppia di punti in $X$ esiste una geodetica che li congiunge.\\
	Sia $r>0$. Lo spazio $X$ è detto \textit{$r$-geodetico} se $\forall x,y\in X$ tali che $d(x,y)<r$, esiste una geodetica che congiunge $x$ e $y$.\\
	Uno spazio geodetico $X$ è \textit{unicamente geodetico} se per ogni coppia di punti, la geodetica che li congiunge (nell'ordine fissato) è unica. In tal caso, l'unico segmento geodetico tra due punti $x,y\in X$ è denotato con il simbolo $[x,y]$. (Per abuso di notazione, un segmento geodetico fissato che congiunge due punti $p$ e $q$ in uno spazio metrico arbitrario sarà ugualmente denotato con il simbolo $[p,q]$.)
\end{defn}
\begin{defn}
	Un sottoinsieme $C\subset X$ è \textit{convesso} se per ogni coppia di punti in $C$ esiste una geodetica che li congiunge, e l'immagine di ogni tale geodetica è contenuta in $C$. 
\end{defn}
\begin{defn}
	Un \textit{triangolo geodetico} in $X$ con vertici $\{p,q,r\}\subset X$ è un sottoinsieme $\Delta(p,q,r)=\vert c_{pq} \vert\cup\vert c_{qr}\vert\cup\vert c_{pr}\vert$, ove $c_{xy}$ è una geodetica che congiunge i punti $x,y\in X$. 
\end{defn}

Denotiamo con $\mathbb E^n$, $n>0$, lo \textit{spazio euclideo}. Come insieme, $\mathbb E^n=\mathbb R^n$ e la metrica $d=d_\E$ è quella euclidea: \[d(x,y)=\sqrt{\sum_{i=1}^n(x_i-y_i)^2},\]
definita per ogni $x=(x_1,\dots,x_n),y=(y_1,\dots,y_n)\in \mathbb E^n$. Lo spazio euclideo è uno spazio metrico unicamente geodetico le cui geodetiche sono le rette. Dati tre numeri reali non-negativi $a,b$ e $c$, esiste un triangolo con lati di lunghezza $a,b$ e $c$ se e solo se $a+b<c$.

\begin{defn}{\label{comparisonEuclide}}
	Sia $\Delta=\Delta(p,q,r)\in X$ un triangolo geodetico in $X$. Un \textit{triangolo euclideo di confronto} per $\Delta$ è un triangolo euclideo $\overline\Delta=\overline{\Delta}(p,q,r)=\Delta(\bar{p},\bar{q},\bar{r})\subset \mathbb E^2$ tale che $d(x,y)=d_{\mathbb E^2}(\bar x,\bar y)$, $\forall x,y\in \{p,q,r\}$. Se $x\in \Delta$ è un punto che giace sul segmento che congiunge $p$ e $q$, un suo \textit{punto di confronto} su $\bar\Delta$ è il punto $\bar x$ che giace sul segmento $[\bar p,\bar q]$ e soddisfa $d(x,p)=d_{\mathbb E^2}(\bar x,\bar p)$. L'angolo interno di $\overline\Delta$ in $\bar p$ è chiamato \textit{angolo di confronto} tra $q$ ed $r$ in $p$ ed è denotato con $\overline\angle_p(q,r)$, oppure con $\angle_p^{(0)}(q,r)$; nel caso in cui $X$ non sia unicamente geodetico, è necessario esplicitare le geodetiche considerate.
\end{defn}
\begin{lem}
	Ogni triangolo geodetico ha un triangolo euclideo di confronto ed esso è unico a meno di isometrie.
	\begin{proof}
		L'unicità è ovvia, poiché un triangolo euclideo è determinato a meno di isometrie dalle lunghezze dei suoi lati. Ora, per la disuguaglianza triangolare in $X$, detto $\Delta=\Delta(p,q,r)$, si ha $d(p,q)\leq d(p,r)+d(r,q)$. Quindi esiste un triangolo euclideo avente i lati di lunghezze $d(p,q),d(p,r),d(r,q)$. 
	\end{proof}
\end{lem}
Esiste un'altra nozione di angolo, dovuta al matematico russo Alexandrov, che consente di parlare di angolo compreso fra due geodetiche. Nel caso di un triangolo geodetico, l'angolo interno in un vertice del triangolo (nel senso della Definizione \ref{comparisonEuclide}) coincide con l'angolo di Alexandrov tra le geodetiche incidenti nel vertice.
\begin{defn}
	Sia $X$ uno spazio metrico e siano $c:[0,a]\to X$ e $c':[0,a']\to X$ due geodetiche con lo stesso punto iniziale $c(0)=c'(0)$. Per $t\in (0,a]$ e $t'\in(0,a']$, consideriamo il triangolo di confronto euclideo $\overline\Delta(c(0),c(t),c'(t'))$ e l'angolo di confronto $\overline\angle_{c(0)}(c(t),c'(t'))$. L'\textit{angolo (di Alexandrov)} tra le due geodetiche $c$ e $c'$ è il numero $\angle_{c,c'}\in [0,\pi]$ definito da: \[\angle_{c,c'}=\angle(c,c')=\limsup_{t,t'\to 0}\overline\angle_{c(0)}(c(t),c'(t'))=\lim_{\varepsilon\to 0}\sup_{0<t,t'<\varepsilon}\overline\angle_{c(0)}(c(t),c'(t'))\]
\end{defn}
Si può dimostrare \cite[Prop. I 1.14 pag. 10]{HB} che vale una disuguaglianza triangolare per gli angoli. Precisamente, se $c,c'$ e $c''$ sono geodetiche in uno spazio metrico $X$ uscenti da uno stesso punto $p\in X$, allora \begin{align}{\label{angleDis}}
\angle_p(c',c'')\leq \angle_p(c,c')+\angle_p(c,c'').
\end{align}
La nozione di angolo consente di associare uno spazio metrico, detto spazio delle direzioni, ad una coppia $(X,p)$, ove $p$ è un punto dello spazio $X$. Lo spazio delle direzioni è un oggetto che contiene alcune informazioni locali dello spazio nell'intorno del punto. 
\begin{defn}
	Sia $X$ uno spazio metrico e sia $p\in X$. Due geodetiche non banali $c,c'$ con lo stesso punto iniziale $p$ \textit{definiscono la stessa direzione} se l'angolo di Alexandrov tra di loro è nullo. La disuguaglianza triangolare per gli angoli (\ref{angleDis}) implica che [$c\sim c'\Leftrightarrow \angle(c,c')=0$] definisce una relazione di equivalenza. Lo \textit{spazio delle direzioni} (o \textit{link}) di $X$ in $p$ è il quoziente dell'insieme delle geodetiche uscenti da $p$ rispetto alla relazione $\sim$, munito della metrica indotta da $\angle$ sulle classi di equivalenza. 
\end{defn}
\section{Spazi $\cat 0$} 
\begin{defn}
	Uno spazio metrico geodetico $(X,d)$ è detto \textit{spazio} $\cat{0}$ se per ogni triangolo geodetico $\Delta\subset X$, vale $\forall p,q\in \Delta$: $$d(p,q)\leq d_{\mathbb E^2}(\bar p,\bar q),$$ ove $\bar p,\bar q$ sono punti di confronto per $p,q$ su un triangolo di confronto euclideo di $\Delta$.
\end{defn}
Dalla definizione precedente segue che gli spazi $\cat 0$ sono unicamente geodetici. Per provarlo abbiamo bisogno del seguente lemma.
\begin{lem}[Convessità della metrica]
	Sia $X$ uno spazio $\cat{0}$. Siano $\gamma,\delta:[0,1]\to X$ due geodetiche riparametrizzate. Allora $\forall t\in[0,1]$ si ha:
	\begin{align}
	d(\gamma(t),\delta(t))\leq(1-t)d(\gamma(0),\delta(0))+td(\gamma(1),\delta(1))
	\end{align}
	\begin{proof}
		Poiché $X$ è geodetico, esiste una geodetica riparametrizzata $\alpha:[0,1]\to X$ da $\gamma(0)$ a $\delta(1)$. Per il Teorema di Talete (valido sul piano euclideo \cite[VI.2]{Eu}), valgono: \begin{align*}
		d_\E(\overline{\gamma(t)},\overline{\alpha(t)})=td_\E(\overline{\gamma(1)},\overline{\alpha(1)})\\
		d_\E(\overline{\alpha(t)},\overline{\delta(t)})=(1-t)d_\E(\overline{\alpha(0)},\overline{\delta(0)})
		\end{align*}
		Applicando la stima $\cat{0}$ al triangolo $\Delta(\gamma(0),\gamma(1),\delta(1))$, otteniamo:
		\begin{align*}
		d(\gamma(t),\delta(t))&\leq d(\gamma(t),\alpha(t))+d(\alpha(t),\delta(t))\leq\\
		&\leq d_\E(\overline{\gamma(t)},\overline{\alpha(t)})+d_\E(\overline{\alpha(t)},\overline{\delta(t)})=\\
		&=td_\E(\overline{\gamma(1)},\overline{\alpha(1)})+(1-t)d_\E(\overline{\alpha(1)},\overline{\delta(1)})=\\
		&=td(\gamma(1),\delta(1))+(1-t)d(\gamma(0),\delta(0))
		\end{align*}
	\end{proof}
\end{lem}
\begin{prop}
	Sia $X$ uno spazio $\cat{0}$. Allora $X$ è unicamente geodetico. 
	\begin{proof}
		Siano $x_0,x_1\in X$ e $\gamma,\delta:[0,1]\to X$ due geodetiche riparametrizzate da $x_0$ a $x_1$, i.e. $\gamma(i)=\delta(i)=x_i,\ i\in\{0,1\}$. Per la convessità della metrica, ricaviamo: $$d(\gamma(t),\delta(t))\leq (1-t)d(x_0,x_0)+td(x_1,x_1)=0,$$ ovvero, $\gamma(t)=\delta(t)$.
	\end{proof}
\end{prop}
\begin{exa}
	Il piano euclideo è chiaramente uno spazio $\cat 0$ perché ogni triangolo geodetico.
\end{exa}
	Sia $\Gamma=(V,E)$ un grafo orientato connesso. Ciò significa che $V$ ed $E$ sono due insiemi, rispettivamente l'insieme dei \textit{vertici} e \textit{archi} di $\Gamma$, e sono date due applicazioni $\de_0,\de_1:E\to V$. Inoltre, $V$ è l'unione delle immagini di $\de_0$ e $\de_1$.\\
	Associamo a $\Gamma$ l'insieme $X=X_\Gamma$ ottenuto prendendo il quoziente di $E\times[0,1]$ rispetto alla relazione di equivalenza generata da \begin{center}$(e,i)\sim (e',i') $ se $\de_i(e)=\de_{i'}(e')$,\end{center} ove $e,e'\in E$ e $i,i'\in\{0,1\}$. Per ogni $e\in E$ è definita una funzione $f_e:[0,1]\to X$ che associa a $t\in [0,1]$ la classe di equivalenza di $(e,t)$ in $X$. Sia $\lambda:E\to(0,+\infty)$ una funzione, che associa a un arco $e\in E$ la sua  \textit{lunghezza} $\lambda(e)$. \\
	Un \textit{cammino lineare a tratti} è una funzione $c:[0,1]\to X$ per cui esiste una partizione $0=t_0<t_1<\dots<t_n=1$ tale che $c\vert_{[t_i,t_{i+1}]}=f_{e_i}\circ c_i$, ove $e_i\in E$ e $c_i$ è una mappa affine $[t_i,t_{i+1}]\to [0,1]$. Diremo che $c$ congiunge $c(0)$ e $c(1)$. La \textit{lunghezza} di $c$ è il numero $l(c):=\sum_{i=0}^{n-1}l(c_i)$, ove $l(c_i)=\lambda(e_i)\vert c_i(t_i)-c_{i+1}(t_{i+1})\vert$.\\
	Possiamo a questo punto definire una pseudometrica $d$ su $X$ ponendo per ogni $x,y\in X$: \begin{align}{\label{grafoMetrico}}d(x,y)=\inf\{l(c)\vert\ c\mbox{\ è\ un\ cammino\ lineare\ a\ tratti\ da}\ x \mbox{\ a\ }y\}.\end{align}
	L'insieme $X$ con la pseudometrica $d$ è detto \textit{grafo metrico} associato al grafo $\Gamma$ e alla funzione di lunghezza $\lambda$.
\begin{exa}
	Sia $\Gamma$ un \textit{albero} metrico, cioè un grafo metrico privo di cicli. Attribuiamo ad ogni arco di $\Gamma$ una lunghezza pari a $1$. L'albero con tale metrica è uno spazio metrico. Inoltre, esso è $\cat 0$. Infatti, per definizione, $\Gamma$ è connesso e la distanza tra due punti è pari alla lunghezza dell'unico cammino privo di inversione che li congiunge. Quindi $\Gamma$ è (unicamente) geodetico. Siano poi $v,v',v''$ tre punti qualsiasi in $\Gamma$. Se un punto è contenuto nella geodetica congiungente gli altri due, il triangolo da essi definito è degenere e soddisfa quindi la disuguaglianza $\cat 0$. Siano quindi $v,v',v''$ non allineati. Il triangolo geodetico che li ha come estremi è omeomorfo ad una stella con tre foglie. Anche in questo caso, dunque, vale la disuguaglianza $\cat 0$.
\end{exa}

\chapter{Teorema del punto fisso di Bruhat-Tits}
In questo capitolo vediamo una prima applicazione degli spazi $\cat 0$. L'esposizione seguente è tratta principalmente da \cite{L}.\\
 Consideriamo un'azione di un gruppo su uno spazio $\cat 0$.  Il Teorema di Bruhat-Tits si occupa del problema di determinare l'esistenza di un punto fisso per tale azione. 
 \section{Spazi di Bruhat-Tits}
\begin{prop}[Disuguaglianza CN di Bruhat-Tits]
	Sia $X$ uno spazio metrico geodetico. Allora $X$ è uno spazio $\cat 0$ se e solo se vale la seguente condizione: \begin{center}{\label{CN}}
		$\forall p,q,r\in X$, $\forall m\in X$ con $d(q,m)=d(r,m)=d(q,r)/2$, si ha \begin{equation}{\label{CNDis}}
		d(p,q)^2+d(p,r)^2\geq 2d(m,p)^2+2d(m,q)^2.\tag{CN}
		\end{equation}
	\end{center} 
	\begin{proof}
		Osserviamo anzitutto che in $\mathbb E^2$ vale l'uguaglianza (\ref{CNDis}). Infatti, siano $P,Q,R,M\in\mathbb E^2$ tali che $d(Q,M)=d(R,M)=d(Q,R)/2$; avendo definito i vettori $x=M-Q, y=x=R-M, u=R-P, z=Q-P,t=M-P$ si ha:
		\begin{align}\begin{split}
		{\label{CNinE2}}
		&z^2+u^2=\vert t-x\vert ^2+\vert t+y\vert^2=\\
		=&t^2+x^2-2(t,x)+t^2+y^2+2(t,y)=\\
		=&2t^2+x^2+y^2=2t^2+2x^2.
		\end{split}\end{align}
		Sia $X$ uno spazio $\cat 0$. Allora, denotando con $\overline z$ il punto di confronto in $\E$ di un generico punto $z$ del triangolo $\Delta(p,q,r)$, si ha: \begin{align*}&d(p,q)^2+d(p,r)^2=d(\bar p,\bar q)^2+d(\bar p,\bar r)^2\overset{(\ref{CNinE2})}{=}\\
		=&2d(\bar m,\bar p)^2+2d(\bar m,\bar q)^2\overset{\cat 0}{\geq}\\
		=&2d(m,p)^2+2d(m,q)^2,\end{align*} poiché $\bar m$ è punto medio del segmento $[\bar q,\bar r]$.\\
		%
		Viceversa, sia $X$ uno spazio geodetico in cui vale \ref{CNDis}. Siano $p,q,r\in X$ i vertici di un triangolo geodetico $\Delta$ in $X$. Sia $x\in [q,r]$ un punto di un lato del triangolo $\Delta$. Possiamo supporre che $d(x,q)\leq d(x,r)$. Sia dunque $r'$ l'unico punto di $[q,r]$ distinto da $q$ a distanza $d(q,x)$ da $x$. Siano ora $\bar p,\bar q,\bar r'$ i punti di confronto in $\E$ rispettivamente per $p,q$ e $r'$. Allora vale 
		\begin{align*}
		d(x,p)^2&\leq \frac1{2}\left[d(p,q)^2+d(p,r')^2-\frac1{2} d(q,x)^2\right]=\\
		&=\frac1{2}\left[d(\bar p,\bar q)^2+d(\bar p,\bar r')^2-\frac1{2} d(\bar q,\bar x)^2\right]=\\
		&=d(\bar x,\bar p)^2
		\end{align*} 
	\end{proof}
\end{prop}
\begin{defn}
	Sia $X$ uno spazio metrico. Diciamo che $X$ soddifsa la \textit{legge del parallelogramma} se $\forall x_1,x_2\in X$, $\exists z\in X$ tale che $\forall x\in X$, \[d(x_1,x_2)^2+4d(z,x)^2\leq 2d(x,x_1)^2+2d(x,x_2)^2.\] 
	Uno spazio completo $X$ che soddisfa tale legge è detto \textit{spazio di Bruhat-Tits}.
\end{defn}
Per $x=x_i$ si ricava $2d(x_i,z)\leq d(x_1,x_2)$; per la disuguaglianza triangolare, $d(x_1,x_2)\leq d(x_1,z)+d(z,x_2)$, dunque, $d(z,x_1)=d(z,x_2)=d(x_1,x_2)/2$. Inoltre, $z$ è univocamente determinato da $x_1,x_2$ poiché se $z'$ è un altro tale punto, ponendo $x=z'$ si ricava $4d(z',x_1)^2+4d(z',z)^2\leq 4d(z',x_1)$. Per questo motivo, $z$ è detto il \textit{punto medio} di $x_1$ e $x_2$.
\begin{cor}
	Sia $X$ uno spazio $\cat 0$ completo.\footnote{Gli spazi $\cat 0$ completi sono detti \textit{spazi di Hadamard}} Allora $X$ è uno spazio di Bruhat-Tits.
	\begin{proof}
		Siano $q,r\in X$ e sia $m$ il punto medio del segmento $[q,r]$. Allora $\forall p\in X$, vale la disuguaglianza (\ref{CNDis}), cioè \[2d(m,p)^2+2d(m,q)^2\leq d(p,q)^2+d(p,r)^2.\] Ora, essendo $2d(m,q)=d(q,r)$, moltiplicando per $2$ ambo i membri della disuguaglianza, ricaviamo:\[d(q,r)^2+4d(m,p)^2\leq 2d(p,q)^2+2d(p,r)^2.\]
	\end{proof}
\end{cor}
\section{Teorema di punto fisso}
Il seguente Teorema è dovuto a J.P. Serre ed è necessario per la dimostrazione del Teorema di Punto Fisso.
\begin{lem}[Teorema di Serre]{\label{lemBru}}
	Sia $X$ uno spazio di Bruhat-Tits. Sia $S\subseteq X$ un sottoinsieme limitato. Allora esiste un'unica palla chiusa di raggio minimo contenente $S$.
	\begin{proof}
		Proviamo anzitutto l'unicità. Supponiamo che $\bar B_r(x_1),\bar B_r(x_2)$ siano due palle di raggio minimo contenenti $S$ e centrate in punti distinti $x_1,x_2\in X$. Sia $x\in S$; vale allora $d(x,x_i)\leq r$, $i=1,2$. Sia $z$ il punto medio di $x_1$ e $x_2$. Per la legge del parallelogramma, \[d(x_1,x_2)^2\leq 2d(x_1,x)^2+2d(x_2,x)^2-4d(x,z)^2\leq 4r^2-4d(x,z)^2. \]
		Poiché $r$ è minimo, $\forall \varepsilon>0$ $\exists x\in S$ tale che $d(x,z)\geq r-\varepsilon$ (altrimenti per qualche $\varepsilon>0$, si avrebbe $S\subseteq \bar B_{r-\varepsilon}(z)$.) Quindi $\forall \varepsilon>0$, si ha $d(x_1,x_2)^2\leq 4r^2-4(r-\varepsilon)^2$, che per $\varepsilon\to 0$ converge a $0$. Segue che $x_1=x_2$ e dunque $\bar B_r(x_1)=\bar B_r(x_2)$.\\
		Proviamo l'esistenza. Sia $r$ l'estremo inferiore dei raggi $r'$ per cui $S$ è contenuto in una palla chiusa di raggio $r'$. Sia $r_n$ una successione reale che converge dall'alto a $r$, i.e. $r_n\overset{n\to\infty}{\longrightarrow}r,\ r_n>r$. Siano $x_n\in X$ i centri delle palle di raggio $r_n$ per cui $S\subseteq \bar B_{r_n}(x_n)$.\\
		La successione $(x_n)$ è di Cauchy. Sia $z_{mn}$ il punto medio di $x_n$ e $x_m$. Per definizione di $r$, $\forall\varepsilon>0\ \exists x\in S$ tale che $$d(x,z_{mn})^2\geq r^2-\varepsilon.$$
		Applicando la legge del parallelogramma con $z=z_{mn}$ ricaviamo \begin{align*}
		d(x_m,x_n)^2\leq& 2d(x,x_m)^2+2d(x,x_n)^2-4d(x,z_{mn})^2\leq\\
		\leq&2r_m^2+2r_n^2-4r^2+4\varepsilon
		\end{align*} Per $m,n\to +\infty$, $2r_m^2+2r_n^2-4r^2\to 0$, ovvero $d(x_m,x_n)\to 0$. \\
		Per completezza di $X$, la successione di Cauchy $x_n$ converge a un punto $x\in X$. Ora, se $s$ è un punto di $S$, vale \begin{align*}d(x,s)\leq d(x,x_n)+d(x_n,s)\leq d(x,x_n)+r_n\overset{n\to\infty}{\longrightarrow} r,\end{align*} cioè $s\in \bar B_r(x)$. 
	\end{proof}
\end{lem}
\begin{defn}
	Il centro della palla di raggio minimo contenente l'insieme limitato $S$ si chiama \textit{circocentro} di $S$.
\end{defn}
\begin{thm}[Teorema del punto fisso di Bruhat-Tits]
	Sia $G$ un gruppo di isometrie per uno spazio di Bruhat-Tits $X$. Se $G$ ha un'orbita limitata, allora $G$ ammette un punto fisso (per esempio il circocentro di questa orbita).
	\begin{proof}
		Sia $S=Gx$ un'orbita limitata dell'azione di $G$ su $X$. Sia $x_0$ il circocentro di $S$ e $B$ la palla chiusa di raggio minimo $r$ contenente $S$. Sia ora $ g\in G$. Allora $gS=S$ e $gB$ è una palla chiusa di centro $x_0'=gx_0$ e raggio uguale a quello di $B$. Per l'unicità asserita nel Lemma \ref{lemBru}., si ha $x_0'=x_0$, e dunque $gx_0=x_0$.
	\end{proof}
\end{thm}

\newpage
\chapter{Spazi CAT($\ka$)}
Nel primo capitolo abbiamo definito la nozione di spazio $\cat 0$ ricorrendo ad un confronto tra i triangoli di uno spazio con quelli dello spazio euclideo. Lo spazio euclideo ha la proprietà di essere una varietà riemanniana piatta su cui il suo gruppo di isometrie agisce transitivamente. Esistono però altri spazi che godono di quest'ultima proprietà di simmetria che possono essere usati come modelli per studiare la natura dei triangoli in uno spazio metrico arbitrario.
\section{Spazi modello a curvatura costante}
Definiamo i seguenti tre spazi metrici fondamentali.
\begin{defn}
	Sia $n>0$ un numero naturale. \begin{enumerate}
		\item $\mathbb E^n$ è lo spazio metrico $\mathbb R^n$ con la metrica euclidea: \[d_{\mathbb E^n}(x,y)=\sqrt{\sum_{i=1}^n(x_i-y_i)^2}.\]
		Essa è indotta dal prodotto scalare euclideo definito ponendo $\forall x,y\in \mathbb E^n$: \[(x\vert y)=\sum_{i=1}^n x_iy_i.\]
		\item $\mathbb S^n$ è la $n$-\textit{sfera unitaria}, ovvero l'insieme \[\{x\in \mathbb R^{n+1}\vert\ (x\vert x)=1\},\] munito della metrica $d_{\mathbb S^n}:\mathbb S^n\times\mathbb S^n\to \mathbb R$ che associa a una coppia $(A,B)$ di punti di $\mathbb S^n$ l'unico numero reale $ d(A,B)\in[0,\pi]$ tale che $$\cos d(A,B)=(A\vert B).$$
		\item $\mathbb H^n$ è l'$n$-\textit{spazio iperbolico}, ovvero l'insieme \[\{x\in \mathbb R^{n+1}\vert\ -x_{n+1}^2+\sum_{i=1}^nx_i^2=-1,\ x_{n+1}>0\}\] munito della metrica $d_{\mathbb H^n}:\mathbb H^n\times\mathbb H^n\to \mathbb R$ che associa a una coppia $(A,B)$ di punti di $\mathbb H^n$ l'unico numero reale $d(A,B)\geq 0$ tale che $$\cosh d(A,B)=A_{n+1}B_{n+1}-\sum_{i=1}^nA_iB_i.$$ Questa metrica è indotta dal prodotto scalare con segnatura $(n,1,0)$ definito su $\mathbb E^{n,1}=\mathbb R^n\times \mathbb R$ ponendo $\forall x,y\in \mathbb E^{n,1}$: $$\langle x\vert y \rangle =-x_{n+1}y_{n+1}+\sum_{i=1}^n x_iy_i.$$ 
	\end{enumerate} 
\end{defn}
\begin{defn}
	Dato un numero reale $\ka$, denotiamo con $\mnk$ i seguenti spazi metrici:
	\begin{enumerate}
		\item se $\ka=0$, $M^n_0$ è lo spazio euclideo $\mathbb E^n$;
		\item se $\ka>0$, $\mnk$ è la $n$-sfera $\mathbb S^n$ con la metrica $\frac{1}{\sqrt\ka}d_{\mathbb S^n}$;
		\item se $\ka<0$, $\mnk$ è lo spazio iperbolico $\mathbb H^n$ con la metrica $\frac{1}{\sqrt{-\ka}}d_{\mathbb H^n}$.
	\end{enumerate}	
	Chiameremo \textit{diametro} di $\mnk$ il numero $D_\ka=\pi/\sqrt{\ka}$ se $\ka>0$ e $D_\ka=\infty$ se $\ka\leq 0$. Il parametro $\ka$ si chiama \textit{curvatura} di $\mnk$ e coincide con la curvatura sezionale (costante) della varietà riemanniana $\mnk$.
\end{defn}
Nel prossimo Teorema richiamiamo alcune proprietà degli spazi modello $\mnk$. Le dimostrazioni non presentate si possono trovare in \cite{HB}.
\begin{thm}{\label{cos}}\begin{itemize}
		\item	$\mnk$ è uno spazio metrico geodetico. 
		\item Se $\ka\leq 0$, $\mnk$ è unicamente geodetico e ogni palla in $\mnk$ è convessa. Se $\ka>0$, esiste un'unica geodetica che congiunge $x,y\in \mnk$ se e solo se $d(x,y)<D_\ka$. Se $\ka>0$, ogni palla chiusa in $\mnk$ di raggio $D_\ka/2$ è convessa.
	\item Inoltre, valgono le seguenti $\operatorname{leggi\ del\ coseno}$: dato un triangolo geodetico in $\mnk$ con lati di lunghezza positiva $a,b$ e $c$ e detto $\gamma$ l'angolo al vertice opposto al lato di lunghezza $c$, si ha:
	\begin{enumerate}
		\item se $\ka=0$: \[c^2=a^2+b^2-2ab\cos(\gamma)\]
		\item se $\ka<0$: 
		\[\cosh(\sqrt{-\ka}c)=\cosh(\sqrt{-\ka}a)\cosh(\sqrt{-\ka}b)-\sinh(\sqrt{-\ka}c)\sinh(\sqrt{-\ka}c)\cos(\gamma)\]
		\item se $\ka>0$: \[\cos(\sqrt{\ka}c)=\cos(\sqrt{\ka}a)\cos(\sqrt{\ka}b)-\sin(\sqrt{\ka}c)\sin(\sqrt{\ka}c)\cos(\gamma).\]
	\end{enumerate}
\end{itemize}
\end{thm}
\begin{lem}[Esistenza di triangoli di confronto in $M_\ka^2$]
	Sia $\ka\in\mathbb R$ e sia $X$ uno spazio metrico. Siano $p,q,r\in X$ tali che $d(p,q)+d(p,r)+d(q,r)<2D_\ka$ (ove $D_\ka=\pi/\sqrt{\ka}$, e se $\ka\leq 0$, la condizione è vuota, ovvero è soddisfatta da ogni terna di punti). Allora esistono dei punti $\bar p,\bar q,\bar r\in M_\ka^2$ tali che \[d(p,q)=d_{M_\ka^2}(\bar p,\bar q),\ d(q,r)=d_{M_\ka^2}(\bar q,\bar r),\ d(r,p)=d_{M_\ka^2}(\bar r,\bar p).\]
	\begin{proof}
		Siano $a=d(p,q)$, $b=d(q,r)$ e $c=d(r,p)$. Possiamo supporre che $a\leq b\leq c$. Per la disuguaglianza triangolare, $c\leq a+b$, quindi $c<D_\ka$ (altrimenti, $a+b+c\geq c+c\geq 2D_\ka$, contro l'ipotesi). Per la legge del coseno (Teorema \ref{cos}.), possiamo trovare un numero $\gamma\in [0,\pi]$ che risolve l'equazione per $M_\ka^2$. Fissato un punto $\bar p\in M_\ka^2$, costruiamo due segmenti geodetici da esso uscenti $[\bar p,\bar q],[\bar p,\bar r]$ di lunghezza rispettivamente $a$ e $b$ che formano un angolo di ampiezza $\gamma$ in $\bar p$. Per la legge del coseno, $d_{M_\ka^2}(\bar q,\bar r)=c$. 
	\end{proof}
	\begin{defn}
		Il triangolo geodetico $\Delta$ di vertici $\bar p,\bar q,\bar r$ si chiama \textit{triangolo di confronto} per la terna $(p,q,r)$, o per un qualsiasi triangolo geodetico $\Delta(p,q,r)$ in $X$. Se $x\in \Delta$ giace sul segmento geodetico $[p,q]$, il punto $\bar x\in\bar\Delta$ si dice \textit{punto di confronto} di $x$ se $\bar x$ giace sul segmento $[\bar p,\bar q]$ e valgono $d(x,p)=d_\mnk (\bar x,\bar p)$,  $d(x,q)=d_\mnk (\bar x,\bar q)$.
	\end{defn}
	\begin{defn}
		Siano $X$ uno spazio metrico geodetico e $\ka\in\mathbb R$. Siano $p,q,r\in X$ tre punti distinti con $d(p,q)+d(p,r)+d(q,r)<2D_\ka$. L'\textit{angolo di $\ka$-confronto} tra $q$ ed $r$ in $p$, denotato $\angle_p^{(\ka)}(q,r)$, è l'angolo in $\bar p$ in un triangolo di confronto $\Delta(\bar p,\bar q,\bar r)\subset M_\ka^2$ per $(p,q,r)$. (In particolare, $\angle_p^{(0)}(q,r)=\overline\angle_p(q,r)$.) 
	\end{defn}
\end{lem}
\section{Spazi CAT$(\kappa)$}
\begin{defn}\footnote{Il termine $\mbox{CAT}$ fu coniato da M. Gromov ed è l'acronimo dei matematici E. Cartan, A.D. Alexandrov e V.A. Toponogov, i quali si occuparono di studiare spazi di questo genere.}
	Sia $X$ uno spazio metrico e sia $\kappa\in\mathbb R$. Sia $\Delta\subset X$ un triangolo geodetico di perimetro minore di $2D_{\ka}$. Sia $\bar\Delta\subset \mnk$ un triangolo di confronto per $\Delta$. Allora, $\Delta$ soddisfa la \textit{disuguaglianza} $\catk$ se $\forall x,y\in \Delta$ e $\forall \bar x,\bar y\in\bar\Delta$ punti di confronto, si ha $$d(x,y)\leq d_\mnk(\bar x,\bar y).$$ 
	Se $\ka\leq 0$, $X$ è uno \textit{spazio} $\catk$ se $X$ è geodetico e tutti i triangoli geodetici in $X$ soddisfano la disuguaglianza $\catk$.\\
	Se $\ka>0$, $X$ è uno \textit{spazio} $\catk$ se $X$ è $D_\ka$-geodetico e ogni triangolo geodetico di perimetro inferiore a $2D_\ka$ soddisfa la disuguaglianza $\catk$.\\
	Si dice che $X$ è uno \textit{spazio a curvatura $\leq \ka$} se è localmente uno spazio $\catk$, i.e. $\forall x\in X$, $\exists r_x>0$ tale che la palla $B_{r_x}(x)$ sia uno spazio $\catk$ (rispetto alla metrica indotta). 
\end{defn}
\begin{prop}{\label{propCatk}}
	Sia $X$ uno spazio $\catk$. \begin{enumerate}
		\item Se $x,y\in X$ hanno distanza $d(x,y)<D_\ka$, esiste un'unica geodetica $[x,y]$ che li congiunge; inoltre, le geodetiche variano in modo continuo rispetto ai suoi estremi. 
		\item Ogni palla in $X$ di raggio inferiore a $D_\ka/2$ è convessa.
		\item Ogni palla in $X$ di raggio inferiore a $D_\ka$ è contraibile. 
		\item $\forall \lambda<D_\ka$, $\varepsilon>0$, esiste un numero reale $\delta=\delta(\ka,\lambda, \varepsilon)$ tale che se $m$ è il punto medio di un segmento geodetico $[x,y]\subset X$, con $d(x,y)\leq \lambda$ e se $$\max\{d(x,m'),d(y,m')\}\leq \frac1{2}d(x,y)+\delta$$ allora $d(m,m')<\varepsilon$.
		\begin{proof}
			1. Siano $p,q\in X$ con $d(p,q)<D_\ka$. Siano $[p,q],[p,q]'$ due segmenti geodetici che congiungono $p$ e $q$. Siano $r\in [p,q],\ r'\in[p,q]'$ tali che $d(p,r)=d(p,r')$. Siano $[p,r]$ e $[r,q]$ due segmenti geodetici la cui concatenazione è $[p,q]$. Ogni triangolo di confronto in $\mnk$ per $\Delta=[p,q]'\cup[p,r]\cup[r,q]$ è degenere, essendo $d(p,q)=d(p,r)+d(r,q)$. Perciò i punti di confronto $\bar r,\bar r'$ risp. per $r$ e $r'$ coincidono in $\mnk$. Per la disuguaglianza $\catk$, vale $d(r,r')\leq d_\mnk(\bar r,\bar r')=0$, ovvero $r=r'$.\\
			\`E un fatto noto ($[HB]$) che negli spazi $M_\ka^2$ le geodetiche parametrizzate uscenti da un punto si distanziano in modo controllato rispetto alla loro lunghezza. Più precisamente, $\forall \ell<D_\ka$, $\exists C=C(\ell,\ka)$ tale che, se $ c,c':[0,1]\to M_\ka^2$ sono geodetiche riparametrizzate con $c(0)=c'(0)$ e $l(c(t)),l(c'(t))\leq\ell$, si ha $\forall t\in [0,1]$: $$d(c(t),c'(t))\leq Cd(c(1),c'(1)).$$ Siano quindi $(p_n)$ e $(q_n)$ due successioni in $X$ convergenti rispettivamente a $p$ e a $q$. Se $d(p,q)\leq\ell<D_\ka$, possiamo supporre che \[ d(p,p_n),d(p,q_n)< \ell<D_\ka.\]  Siano  $c,c_n,c'_n:[0,1]\to X$ le geodetiche riparametrizzate di traccia (orientata) rispettivamente $[p,q],[p_n,q_n]$ e $[p,q_n]$. Allora, per la disuguaglianza $\catk$ applicata ai triangoli $\Delta(p,q,p_n)$ e $\Delta(p,q,q_n)$, vale \begin{align*}d(c(t),c_n(t))\leq d(c(t),c'_n(t))+d(c'_n(t),c_n(t))\overset{{\catk}}{=}\\
			C (d(c(1),c'_n(1))+d(c'_n(0),c_n(0)))=C(d(q,q_n)+d(p,p_n))\end{align*}
			Dunque, la geodetica $c_n$ congiungente $p_n$ e $q_n$ converge uniformemente a quella congiuntente $p $ e $q$.\\
			Il punto $2$ segue dal fatto che negli spazi $M_\ka^2$ ogni palla di raggio inferiore a $D_\ka/2$ è convessa (Teorema \ref{cos}). Sia $B\subseteq X$ è una palla di raggio $r<D_\ka/2$ centrata in $x_0$. Per ogni $ x,y\in B$, detta $c$ la geodetica che li congiunge, vale \[d(c(t),x_0)\leq d(\overline{c(t)},\overline{x_0})<r\] per la disuguaglianza $\cat\ka$ applicata al triangolo $\Delta(x_0,x,y)$.\\
			3. segue immediatamente da 1. e 2. Se $B=\bar B_\ell (x)$ è una palla chiusa di raggio $\ell<D_\ka$, la mappa $B\times[0,1]\to X$ che associa a $(y,t)$ il punto a distanza $td(x,y)$ da $y$ lungo la geodetica $[x,y]$ è una ritrazione continua di $B$ in $x$.   \\
			4. \cite[Prop. II 1.4(5), pag.161]{HB}.
		\end{proof}
	\end{enumerate}
\end{prop}
La disuguaglianza $\catk$ può essere riformulata in molti modi. La seguente proposizione caratterizza gli spazi $\catk$.
\begin{prop}{\label{1.7}}
	Sia $\ka\in\mathbb R$. Sia $X$ uno spazio metrico $D_\ka$-geodetico. Allora le seguenti condizioni sono equivalenti (quando $\ka>0$, supponiamo che il perimetro di ogni triangolo considerato sia inferiore a $2D_\ka$):
	\begin{enumerate}
		\item $X$ è uno spazio $\catk$.
		\item Per ogni triangolo geodetico $\Delta([p,q],[q,r],[r,p])$ in $X$ e ogni punto $x\in [q,r]$, la seguente disuguaglianza è soddisfatta da ogni punto di confronto $\bar x\in [\bar q,\bar r]\subset\Delta(\bar p,\bar q,\bar r)\subset M_\ka^2$:\[d(p,x)\leq d_\mnk(\bar p,\bar x)\]
		\item Per ogni triangolo geodetico $\Delta([p,q],[q,r],[r,p])$ in $X$ e ogni coppia di punti $x\in [p,q]\setminus\{p\}$, $y\in [p,r]\setminus \{p\}$, gli angoli al vertice corrispondenti a $p$ in $\overline\Delta(p,q,r)\subset M_\ka^2$ e in $\overline{\Delta}(p,x,y)\subset M_\ka^2$ soddisfano: \[\angle_p^{(\ka)}(x,y)\leq\angle_p^{(\ka)}(q,r).\]
		\item L'angolo di Alexandrov tra due lati di un qualunque triangolo geodetico in $X$ con vertici distinti non è più grande dell'angolo tra i corrispondenti lati del suo triangolo di confronto in $M_\ka^2$. 
	\end{enumerate}
\begin{proof}
	\cite[Prop. II 1.7]{HB}.
\end{proof}
%
	
\end{prop}

\begin{cor}
	Per $\ka\leq 0$, ogni spazio $\catk$ è contraibile; in particolare, è semplicemente connesso.
	\begin{proof}
		Sia $x_0$ un punto fissato di $X$. Denotiamo con $c_x:[0,1]\to X$ l'unica geodetica a velocità costante che congiunge $x$ a $x_0$. Per la Proposizione \ref{propCatk}(1), $c{\_}$ varia in modo continuo con gli estremi; perciò la funzione $H:X\times[0,1]\to X$ definita ponendo $\forall x\in X,t\in[0,1]$ \[H(x,t)=c_x(t)\] è continua. Inoltre, $H(x,0) =x$ e $H(x,1)=x_0$, $\forall x\in X$.
	\end{proof}
\end{cor}
\begin{thm}
	Sia $X$ uno spazio metrico. \begin{enumerate}
		\item Se $X$ è $\catk$, allora $X$ è $\cat{\ka'}$ $\forall \ka'\geq\ka$.
		\item Se $X$ è $\cat{\ka'},$ $\forall\ka'>\ka$, allora $X$ è $\catk$. 
	\end{enumerate}
	\begin{proof}Discende immediatamente dal fatto che i triangoli in $\mnk$ soddisfano la disuguaglianza $\cat{\ka'}$ per tutti i $\ka'\geq \ka$ (\cite{HB}). 
	\end{proof}
	
\end{thm}

\section{Coni: il Teorema di Berestovskii}
Lo scopo di questo paragrafo è la dimostrazione del Teorema di Berestovskii  che fornisce un  legame tra gli spazi $\cat 1$ e gli spazi $\catk$. Il Teorema riveste un ruolo fondamentale nella dimostrazione del Criterio di Gromov.\\
Innanzitutto diamo la definizione di cono su uno spazio metrico.
\begin{defn}
	Siano $Y$ uno spazio metrico e $\ka\in\mathbb R$. Il $\ka$-\textit{cono} $X=C_\ka Y$ è definito come segue.\\
	Se $\ka\leq 0$, l'insieme $X$ è il quoziente di $[0,\infty)\times Y$ rispetto alla relazione di equivalenza generata da: $(t,y)\sim (t',y')$ se $t=t'=0$.\\
	Se $\ka>0$, $X$ è il quoziente di $[0,D_\ka/2]\times Y$ rispetto alla stessa relazione di equivalenza.\\
	La classe di $(t,y)$ è denotata con $ty$ e la classe di $(0,y)$ con $0$, ed è chiamata \textit{vertice del cono}.\\	
	\indent Sia $d_\pi(y,y'):=\min\{\pi,d(y,y')\}$. La distanza tra due punti del cono $x=ty$ e $x'=t'y'$ è definita come segue.\\
	Se $\ka=0$, $d(x,x')$ è l'unico numero reale non-negativo tale che \begin{align}d(x,x')^2=t^2+t'^2-2tt'\cos(d_\pi(y,y')).\end{align}
	Se $\ka<0$, $d(x,x')$ è l'unico numero reale non-negativo tale che \begin{align}\begin{split}\cosh(d(x,x')\sqrt{-\ka})=\cosh(t\sqrt{-\ka})&\cosh(t'\sqrt{-\ka})-\\
	&\sinh(t\sqrt{-\ka})\sinh(t'\sqrt{-\ka})\cos(d_\pi(y,y')).\end{split}\end{align}
	Se $\ka>0$, la distanza tra $x$ e $x'$ è l'unico numero $d(x,x')\in [0,D_\ka]$ tale che \begin{align} \cos(d(x,x')\sqrt{\ka})=\cos(t\sqrt{\ka})\cos(t'\sqrt{\ka})-\sin(t\sqrt{\ka})\sin(t'\sqrt{\ka})\cos(d_\pi(y,y')).\end{align}
	 
\end{defn}
In \cite[Prop. I 5.9]{HB}, gli autori dimostrano che le formule precedenti definiscono una metrica sul $\ka$-cono e che $X=C_\ka Y$ è completo se e solo se lo è $Y$. Inoltre, la distanza tra due punti $x=ty$ e $x'=t'y'$ in $X$ è definita in modo tale che $d(x,x')=t$ se $x'=0$, e $\angle_0^{(\ka)}(x,x')=d_\pi(y,y')$ se $t,t'>0$. \\
I due risultati seguenti, che presentiamo senza dimostrazione, mostrano il comportamento dell'applicazione che associa ad uno spazio metrico il suo \textit{$\ka$-cono}. Il primo motiva la definizione di $\ka$-cono, mentre il secondo caratterizza le geodetiche del cono. 
\begin{prop}{\label{coneSphere}}\cite[Prop. I 5.8]{HB}
	Sia $Y=\mathbb S^{n-1}$ una sfera. Allora $X=C_\ka Y$ è isometrico a $\mnk$ se $\ka\leq 0$ e a una palla chiusa di raggio $D_\ka/2$ in $\mnk$ se $\ka>0$. 
\end{prop}
\begin{lem}[Caratterizzazione delle geodetiche]{\label{I.5.10}}\cite[Prop. I 5.10]{HB}  Siano $x_1=t_1y_1,\ x_2=t_2y_2\in C_\ka Y$.
	\begin{enumerate}
		\item Se $t_1,t_2>0$ e $d(y_1,y_2)<\pi$, esiste una bijezione tra l'insieme dei segmenti geodetici congiungenti $y_1$ a $y_2$ in $Y$ e l'insieme dei segmenti geodetici che congiungono $x_1$ a $x_2$ in $C_\ka Y$.
		\item Negli altri casi, esiste una geodetica congiungente $x_1$ e $x_2$; questo segmento è unico, eccetto nel caso in cui $\ka>0$ e $d(x_1,x_2)=D_\ka$.
		\item Ogni segmento geodetico che congiunge $x_1$ a $x_2$ è contenuto in una palla chiusa di raggio $\max\{t_1,t_2\}$ centrata nel vertice $0\in C_\ka Y$. 
	\end{enumerate} 
\end{lem}
\begin{thm}[Berestovskii]{\label{Ber}}
	Sia $Y$ uno spazio metrico. Il $\ka$-cono $X=C_\ka Y$ su $Y$ è $\catk$ se e solo se $Y$ è uno spazio $\operatorname{CAT}(1)$.
	\begin{proof}
		Supponiamo che $Y$ sia uno spazio $\cat 1$ e proviamo che il cono $X$ è uno spazio $\cat \ka$. Siano $x_i=t_iy_i$ tre punti di $X$, tali che il triangolo geodetico $\Delta_X$ da loro definito abbia perimetro inferiore a $2D_\ka$. Se uno dei $t_i$ è $0$, ovvero se $\Delta_X$ ha un vertice nel vertice del cono, il triangolo $\Delta_X$ con la metrica indotta è isometrico al suo triangolo di confronto in $M_\ka^2$.\footnote{Le geodetiche $[0,x_i]$ hanno la forma $t\in[0,t_i]\mapsto ty_i$ e, per esempio nel caso euclideo ($\ka=0$), $d(ty_1,t'y_2)=t^2+t'^2-2tt'\cos\angle_0^{(0)}(x_1,x_2)=t^2+t'^2-2tt'\cos\angle_0(\bar x_1,\bar x_2)=d(\overline{ty_1},\overline{t'y_2})$ per la legge del coseno.  }
		Supponiamo dunque $t_i>0$, $\forall i$.
		Consideriamo tre casi:
		\begin{enumerate}
			\item $d(y_1,y_2)+d(y_2,y_3)+d(y_3,y_1)<2\pi$,
			\item $d(y_1,y_2)+d(y_2,y_3)+d(y_3,y_1)\geq 2\pi$, ma $d(y_i,y_j)<\pi$ $\forall i,j=1,2,3$,
			\item per qualche $i,j$, $d(y_i,y_j)\geq \pi$.
		\end{enumerate} 
		Caso 1. Sia $\Delta=\Delta(y_1,y_2,y_3)=[y_1,y_2]\cup [y_2,y_3]\cup[y_3,y_1]\subset Y$.\\ 
		Fissiamo un triangolo di confronto $\overline \Delta$ in $M_1^2=\mathbb S^2$ di vertici $\bar y_1,\bar y_2,\bar y_3$. La mappa di confronto $\overline{\Delta}\to\Delta$, per la Proposizione \ref{coneSphere}, si estende ad una bijezione $C_\ka\overline\Delta\subset M_\ka^3\to C_\ka\Delta$. Sia $x=ty\in[x_2,x_3]$ e sia $\bar y\in [\bar y_2,\bar y_3]$ il punto di confronto per $y$. Posto $\bar x_i =t\bar y_i$ si ha $d_{C_\ka Y}(x_i,x_j)=d_{C_\ka \mathbb M_1^2}(t_i\bar y_i,t_j\bar y_j)=d_{M_k^3}(\bar x_i,\bar x_j)$; dunque il triangolo $\Delta(\bar x_1,\bar x_2,\bar x_3)\subset M_\ka^3$, è un triangolo di $\ka$-confronto per $\Delta_X$, con $\bar x=t\bar y$ come punto di confronto per $x$.\\ Dato che la disuguaglianza $\cat{1}$ vale per $\overline{\Delta}$, abbiamo $d(y_1,y)\leq d(\bar y_1,\bar y)$. Utilizzando la definizione di metrica su $C_\ka Y$, vediamo che $d(x_1,x)\leq d(\bar x_1,\bar x)$. Infatti, se $d(y_1,y)\leq d(\bar y_1,\bar y)$, vale $\cos d_\pi(y_1,y)\geq \cos d_\pi(\bar y_1,\bar y)$ e dunque, per esempio utilizzando la definizione di metrica per il cono euclideo $C_0Y$, abbiamo 
		\begin{align*}
		d(x_1,x)=&t_1^2+t^2-2t_1t\cos d_\pi(y_1,y)\leq\\
		\leq&t_1^2+t^2-2t_1t\cos d_\pi(\bar y_1,\bar y)=\\
		=&d(\bar x_1,\bar x).
		\end{align*}
		Nel caso non-euclideo, la disuguaglianza vale allo stesso modo perché dipende solo dall'addendo che contiene $\cos d_\pi (y_1,y)$. \\
		Caso 2. Siano $\Delta(\tilde{0},\tilde{x}_1,\tilde{x}_2)$ e $\Delta(\tilde{0},\tilde{x}_1,\tilde{x}_3)$ i triangoli di confronto in $M_\ka^2$ rispettivamente per $\Delta(0,{x}_1,{x}_2)$ e $\Delta({0},{x}_1,{x}_3)$, scelti in modo che $\tilde x_2$ e $\tilde{x}_3$ giacciano su lati differenti di $[\tilde{x}_1,\tilde{0}]$. Dalla definizione di metrica su $C_\ka Y$, abbiamo 
		$\angle_{\tilde{0}}(\tilde{x}_1,\tilde{x}_2)=d(y_1,y_2)$, $\angle_{\tilde{0}}(\tilde{x}_1,\tilde{x}_3)=d(y_1,y_3)$,  $\angle_{\tilde{x}_1}(\tilde{0},\tilde{x}_2)=\angle_{x_1}(0,x_2)$ e $\angle_{\tilde{x}_1}(\tilde{0},\tilde{x}_3)=\angle_{x_1}(0,x_3)$.
		 Poiché per ipotesi $d(y_1,y_2)+d(y_1,y_3)>\pi$, \[\angle_{\tilde{0}}(\tilde{x}_2,\tilde{x}_3)=2\pi-\angle_{\tilde{0}}(\tilde{x}_2,\tilde{x}_1)-\angle_{\tilde{0}}(\tilde{x}_3,\tilde{x}_1)\leq d(y_2,y_3)=\angle_0(x_2,x_3) \] e $d(\tilde{x}_2,\tilde{x}_3)\leq d(x_2,x_3)$.
		  Quindi in un triangolo di confronto $\overline\Delta=\Delta(\bar x_1,\bar x_2,\bar x_3)\subset M_\ka^2$ per $\Delta_X=\Delta(x_1,x_2,x_3)$, abbiamo 
		\begin{align*}\angle_{\bar x_1}(\bar x_2,\bar x_3)\geq \angle_{\tilde{x}_1}(\tilde x_2,\tilde x_3)=&\angle_{\tilde{x}_1}(\tilde x_2,\tilde 0)+\angle_{\tilde{x}_1}(\tilde 0,\tilde x_3)=\\=&\angle_{x_1}(x_2,0)+\angle_{x_1}(0,x_3)\geq \angle_{x_1}(x_2,x_3).\end{align*} 
		Applicando la Proposizione \ref{1.7}, segue la tesi.\\
		Caso 3. Supponiamo che $d(y_1,y_3)\geq \pi$. Allora il segmento geodetico $[x_1,x_3]$ è la concatenazione di $[x_1,0]$ e $[0,x_3]$. Siano $\overline{\Delta}_1=\Delta(\bar 0,\bar x_1,\bar x_2)$ e $\overline{\Delta}_3=\Delta(\bar 0,\bar x_3,\bar x_2)$ triangoli di confronto in $M_\ka^2$ rispettivamente per $\Delta_1=\Delta(0,x_1,x_2)$ e $\Delta_3=\Delta(0,x_3,x_2)$, scelti in modo che i vertici $\bar x_1$ e $\bar{x}_3$ giacciano su lati differenti del segmento comune $[0,x_2]$: nel caso $\ka>0$, tali triangoli di confronto esistono perché il perimetro di $\Delta(x_1,x_2,x_3)$ è minore di $2D_\ka$ per ipotesi. La somma degli angoli nel vertice $\bar 0$  in $\overline \Delta_1$ e $\overline{\Delta}_3$ è $d_\pi(y_1,y_2)+d_\pi(y_2,y_3)\geq \pi$. Per il Lemma di Alexandrov \cite[Lemma 2.16, pag.25]{HB} $\Delta(x_1,x_2,x_3)$ soddisfa la disuguaglianza $\catk$.\\
		\indent Rimane da provare che se $X$ è $\catk$, $Y$ è $\cat 1$. Sappiamo dalla Proposizione \ref{I.5.10} che punti a distanza minore di $\pi $ in $Y$ sono congiunti da un'unica geodetica. 
		Siano $\Delta(y_1,y_2,y_3)$ un triangolo in $Y$ di perimetro minore di $2\pi$ e $\overline\Delta=\Delta(\bar y_1,\bar y_2,\bar y_3)$ un triangolo di confronto in $M_1^2=\mathbb S^2$. Dato $y\in [y_2,y_3]$, sia $\bar y\in [\bar y_2,\bar y_3]$ il suo punto di confronto. Sul sotto-cono $C_\ka\Delta$, consideriamo tre punti $x_i=\varepsilon y_i$, per $i=1,2,3$, ove $\varepsilon>0$ è sufficientemente piccola da assicurare che il perimetro del triangolo con vertici $x_1,x_2,x_3$ sia inferiore a $2D_\ka$. Il cono $C_\ka\overline{\Delta}$ è un sotto-cono di $C_\ka M_1^2\subseteq M_\ka^3$ e i punti $\bar x_i=\varepsilon \bar y_i$, $i=1,2,3$, sono i vertici di un triangolo di confronto $\overline {\Delta'}$ per il triangolo geodetico $\Delta'$ con vertici $x_1,x_2,x_3$. Se $x=ty\in[x_2,x_3]$, il suo punto di confronto è $\bar x=t\bar y$. Per ipotesi $d(x_1,x)\leq d(\bar x_1,\bar x)$, perché $X$ è $\catk$. Quindi $d(y_1,y)\leq d(\bar y_1,\bar y)$, per definizione di metrica sul cono $C_\ka Y$.
	\end{proof}
\end{thm}

\newpage
\part{}

\chapter{Complessi}
In questa sezione introduciamo una classe di spazi metrici che, al pari dei complessi di celle, hanno una natura combinatorica. Nella prima parte consideriamo i complessi poliedrici modellati sugli spazi omogenei a curvatura costante. I risultati tecnici in questo ambito che sono necessari per il seguito vengono presentati senza dimostrazione. Rimandiamo il lettore alla consultazione di \cite{HB} per una trattazione completa dell'argomento. \\ Nella seconda parte di questa sezione presentiamo una sotto-classe dei complessi poliedrici che generalizza la nozione di complesso simpliciale euclideo. Il risultato principale è il Teorema \ref{5.18}, che lega la curvatura di un complesso sferico alla struttura dei simplessi che lo compongono.
\section{Complessi poliedrici} 
\begin{defn}
	Sia $\ka\in\mathbb R$ fissato. Una \textit{cella $M_\ka$-poliedrica convessa} $C\subset \mnk$ è l'inviluppo convesso di un insieme finito di punti $P\subset \mnk$; se $\ka>0$, richiediamo che $P$ (e quindi $C$) sia contenuto in una palla aperta di raggio $D_\ka/2$. La \textit{dimensione} di $C$ è la dimensione del più piccolo $m$-piano che lo contiene. L'\textit{interno} di $C$ è la parte interna di $C$ come sottoinsieme di tale $m$-piano. \\
	Sia $H$ un iperpiano in $\mnk$. Se $C$ giace in uno dei due semi-spazi determinati da $H$, e $F:=C\cap H\neq \emptyset$, diremo che $F$ è una \textit{faccia} di $C$; se $F\neq C$, la faccia è detta \textit{propria}. La dimensione della faccia $F$ è la dimensione del più piccolo $m$-piano che la contiene. L'interno di $F$ è l'interno di $F$ come sottoinsieme di tale $m$-piano. Le facce $0$-dimensionali di $C$ sono chiamate \textit{vertici}. Il supporto di un punto $x\in C$ è l'unica faccia di $C$ che contiene $x$ nel suo interno, ed è denotato con $\supp(x)$. 
\end{defn}

\begin{defn}{\label{MkPoli}}
	Sia $(C_\lambda:\ \lambda\in \Lambda)$ una famiglia di celle $M_\ka$-poliedriche convesse $C_\lambda\subset \mm{n_\lambda}$. Sia \[X=\coprod_{\lambda\in\Lambda}C_\lambda=\bigcup_{\lambda\in\Lambda}C_\lambda\times\{\lambda\} \]l'unione disgiunta di tali celle e sia $\sim$ una relazione di equivalenza su $X$. Siano $K=X/\sim$ l'insieme quoziente e $p:X\to K$ la projezione canonica. Sia \[p_\lambda:x\in S_\lambda\mapsto p(x,\lambda)\in K.\] $K$ è un $M_\ka$-\textit{complesso poliedrico} se:
	\begin{enumerate}
		\item $\forall \lambda\in \Lambda$, la restrizione di $p_\lambda$ all'interno di ogni faccia di $C_\lambda$ è injettiva;
		\item $\forall \lambda_1,\lambda_2\in \Lambda$ e $x_i\in C_{\lambda_i}$, se $p_{\lambda_1}(x_1)=p_{\lambda_2}(x_2)$, esiste una isometria $h:\supp(x_1)\to \supp(x_2)$ tale che $p_{\lambda_1}(y)=p_{\lambda_2}(h(y))$, $\forall y\in \supp(x_1)$.\\
		L'insieme delle classi di isometria delle facce delle celle $C_\lambda$ verrà denotato con $\sh K$.
	\end{enumerate} 
\end{defn}
\begin{defn}
	Sia $K$ un complesso $M_\ka$-poliedrico. L'immagine $C=p_\lambda(F)$ di una faccia $n$-dimensionale $F\subset C_\lambda$ prende il nome di $n$-\textit{cella} di $K$; l'interno di $C$ è l'immagine tramite $p_\lambda$ dell'interno della faccia $F$.
	
\end{defn}
Per definire la cosiddetta \textit{pseudometrica intrinseca} di $K$ abbiamo bisogno di una definizione preliminare.
\begin{defn}
	Un \textit{cammino geodetico a tratti} è una mappa $c:[a,b]\to K$ tale che esista una partizione $a=t_0\leq t_1\leq\dots \leq t_k=b$ e geodetiche $c_i:[t_{i-1},t_i]\to C_{\lambda_i}$ tali che $\forall t\in [t_{i-1},t_i]$ si abbia $c(t)=p_{\lambda_i}(c_i(t))$. La \textit{lunghezza} $l(c)$ di $c$ è \[l(c):=\sum_{i=1}^kl(c_i).\]
\end{defn}
\begin{defn}
	Siano $x,y\in K$. La loro distanza $d(x,y)$ è definita come l'estremo inferiore delle lunghezze dei cammini geodetici a tratti che li congiungono. La funzione $d$ è detta \textit{pseudometrica intrinseca} di $K$.
\end{defn}
Per ogni $x\in K$, definiamo il seguente numero:\[\varepsilon(x):=\inf\{\varepsilon(x,C)\vert\ C\subseteq K\mbox{ è una cella contenente }x\}\]ove \[\varepsilon(x,C):=\inf\{d_C(x,F)\vert\ F\mbox{ è una faccia di } C,\ x\notin F\}.\]
\begin{lem}[Bridson]
	Sia $K$ un complesso $M_\ka$-poliedrico. Se $\varepsilon(x)>0$ $\forall x\in K$, allora la pseudometrica intrinseca è una metrica.\\
	Se $\sh{K}$ finito. Allora $\varepsilon(x)>0,$ $\forall x\in K$, e dunque la pseudometrica intrinseca è una metrica. Inoltre, $K$ è uno spazio geodetico completo.
\end{lem}
\begin{lem}{\label{I.7.56}}
	Sia $K$ un complesso $M_\ka$-poliedrico con $\sh K$ finito. Allora ogni punto di $K$ ha un intorno isometrico all'intorno di un qualche punto sufficientemente vicino a un vertice del suo supporto.
\end{lem}
Vedremo nel seguito come la struttura metrica di un complesso $M_\ka$-poliedrico risieda nella natura locale dei suoi vertici. Per descrivere la nozione di \textit{località} abbiamo bisogno di definire alcuni enti geometrici associati ad ogni punto del complesso: il link e la stella (o star).
\begin{defn}
	Sia $x\in K$. Il \textit{link geometrico} $\Lk(x,K)$ è lo spazio delle direzioni in $x$. La \textit{stella aperta} (o \textit{star}) di $x$ in $K$, denotata con $\st(x)$, è l'unione degli interni delle celle di $K$ contenenti $x$. 
\end{defn}

\begin{lem}{\label{7.39}}
	Sia $K$ un complesso $M_\ka$-poliedrico e sia $x\in K$. Se $\varepsilon(x)>0$, $B_{\varepsilon(x)/2}(x)$ è isometrica alla palla aperta di raggio $\varepsilon(x)/2$ centrata nel vertice di $C_\ka(\Lk(x,K))$.
	\begin{proof}
	\cite[I 7.39 pag.115]{HB}
	\end{proof}
\end{lem} 

\begin{thm}{\label{5.4}}
	Sia $K$ un complesso $M_\ka$-poliedrico con $\sh{K}$ finito. \\
	Se $\ka\leq 0$ le seguenti affermazioni sono equivalenti:
	\begin{enumerate}
		\item $K$ è uno spazio $\catk$;
		\item $K$ è unicamente geodetico;
		\item $K$ soddisfa la condizione dei link\footnote{Cfr. Definizione \ref{linkCn}} e non contiene cerchi immersi isometricamente;
		\item $K$ è semplicemente connesso e soddisfa la condizione dei link. 
	\end{enumerate}
	Se $\ka>0$ le seguenti affermazioni sono equivalenti:
	\begin{enumerate}
		\item $K$ è un spazio $\catk$;
		\item $K$ è $D_\ka$-unicamente geodetico;
		\item $K$ soddisfa la condizione dei link e non contiene cerchi di lunghezza inferiore a $2D_\ka$ isometricamente immersi.
	\end{enumerate}
	
\end{thm}
Sia $K$ un complesso euclideo (ovvero $M_0$-poliedrico) le cui celle $C_\lambda$ sono tutte dei cubi $[0,1]^{n_\lambda}$. Un tale complesso è detto \textit{complesso cubico}. Se $K$ soddisfa le seguenti ulteriori assiomi, viene detto \textit{complesso cubico stretto}: \begin{enumerate}
	\item ogni cella $C_\lambda$ è isometrica ad un cubo $[0,1]^{n_\lambda}$,
	\item ogni projezione $p_\lambda$ è injettiva,
	\item l'intersezione di due celle è vuota o è una singola faccia di $K$. 
\end{enumerate}
Per esempio, l'$n$-toro (ottenuto quozientando il cubo $[0,1]^n$ rispetto alla relazione che associa ad ogni punto nelle facce $ (n-1)$-dimensionali la sua projezione ortogonale sulla faccia opposta) è un complesso cubico ma non in senso stretto. Questa classe di complessi avrà un ruolo importante nella prossima sezione: il Teorema di Gromov consente di caratterizzare i complessi cubici che ammettono una metrica $\cat 0$. 
\section{Complessi $M_\ka$-simpliciali}
Come preannuciato, in questa seconda parte del capitolo consideriamo una classe di complessi poliedrici. La definizione generalizza quella di complesso simpliciale euclideo nel caso in cui i simplessi siano sottospazi dello spazio modello $\mnk$ da cui ereditano la metrica. 
\begin{defn}
	Siano $n\leq m$ interi positivi. Un $n$-\textit{piano} in $\mm m$ è un suo sottospazio isometrico a $\mnk$. Diremo che $(n+1)$ punti in $\mm m$ sono in \textit{posizione generale} se non sono contenuti in alcun $(n-1)$-piano.\\
	Un $n$-\textit{simplesso geodetico} $S\subset \mm m$ è l'inviluppo convesso di $(n+1)$ punti in posizione generale, i.e. il più piccolo sottoinsieme convesso di $\mm m$ che contiene tali punti; tali punti sono i \textit{vertici} di $S$. Se $\ka>0$, si richiede che i vertici giacciano in un emisfero aperto, i.e. in una palla aperta di raggio $D_\ka/2$.\\
	Una \textit{faccia} di $S$ è un sottoinsieme $T\subseteq S$ che è l'inviluppo convesso di vertici di $S$; se $T\neq S$, la faccia è detta \textit{propria}. L'\textit{interno} di $S$ è l'insieme dei punti di $S$ che non giacciono in alcuna sua faccia propria.
\end{defn}
\begin{defn}{\label{MkComp}}
	Sia $(S_\lambda:\ \lambda\in \Lambda)$ una famiglia di simplessi geodetici $S_\lambda\subset \mm{n_\lambda}$. Sia \[X=\coprod_{\lambda\in\Lambda}S_\lambda=\bigcup_{\lambda\in\Lambda}S_\lambda\times\{\lambda\} \]l'unione disgiunta di tali simplessi e sia $\sim$ una relazione di equivalenza su $X$. Siano $K=X/\sim$ l'insieme quoziente e $p:X\to K$ la projezione canonica. Sia $p_\lambda:x\in S_\lambda\mapsto p(x,\lambda)\in K$ la restrizione della projezione canonica. $K$ è un $M_\ka$-\textit{complesso simpliciale} se:
	\begin{enumerate}
		\item $\forall \lambda\in \Lambda$, $p_\lambda$ è injettiva;
		\item se $p_\lambda(S_\lambda)\cap p_{\lambda'}(S_{\lambda'})\neq \emptyset$, esistono due facce $T_{\lambda}\subseteq S_{\lambda}$, $T_{\lambda'}\subseteq S_{\lambda'}$ e una isometria $h_{\lambda,\lambda'}:T_\lambda\to T_{\lambda'}$ tali che $p(x,\lambda)=p(x',\lambda')$ se e solo se $x'=h_{\lambda,\lambda'}(x)$. \end{enumerate}
	L'insieme delle classi di isometria delle facce dei simplessi geodetici $S_\lambda$ verrà denotato con $\sh K$.
\end{defn}
\`E chiaro che i complessi $M_\ka$-simpliciali sono complessi $M_\ka$-poliedrici. Le definizioni date nella prima parte sono quindi valide anche in questo caso; ripetiamo le definizioni in questo caso particolare per fissare i termini. 
\begin{defn}
	Sia $K$ un $M_\ka$-complesso simpliciale. L'immagine $S=p_\lambda(S_\lambda)$ di un $m$-simplesso geodetico $S_\lambda$ si chiama $m$-\textit{simplesso} di $K$.\\
	Se $T\subseteq S$ e $p_\lambda^{-1}(T)$ è una faccia del simplesso $p_\lambda^{-1}(S)$, $T$ si chiama \textit{faccia} di $S$. \\
	L'\textit{interno} di $S$ è l'immagine tramite $p_\lambda$ dell'interno di $p_\lambda^{-1}(S)$. \\
	Il \textit{supporto} di un punto $x\in K$ è l'unico simplesso $\supp(x)$ che contiene $x$ nel suo interno.	
\end{defn}

Diamo ora una \textit{definizione simpliciale} di link.
\begin{defn}
	Sia $K$ un complesso $M_\ka$-simpliciale (oppure un complesso simpliciale astratto) e $x$ un suo vertice. I simplessi contenenti $x$ costituiscono un complesso simpliciale astratto i cui vertici sono gli $1$-simplessi contenenti $x$ e, più in generale, gli $(m+1)$-simplessi di $K$ sono gli $m$-simplessi di $K$ contenenti $x$. Tale complesso si indica con $\Lk(x,K)$. Esso può anche essere identificato con il sottocomplesso di $K$ costituito dalle facce  $T$ dei simplessi di $K$ che contengono $x$, ma $x\notin T$.
\end{defn}
\begin{defn}
	Sia $x\in K$. La \textit{stella chiusa}
	di $x$, denotata con $\St(x)$, è l'unione dei simplessi di $K$ che contengono $x$. La \textit{stella aperta}, denotata con $\st(x)$, è l'unione degli interni dei simplessi di $K$ che contengono $x$. \\
	L'immagine tramite $p_\lambda$ di un segmento geodetico $[x_\lambda,y_\lambda]\subseteq T_\lambda\subseteq S_\lambda$ è chiamato \textit{segmento} nel simplesso $T=p_\lambda(T_\lambda)$ ed è denotato con $[x,y]$, ove $x=p_\lambda(x_\lambda),\ y=p_\lambda(y_\lambda)$. La lunghezza di $[x,y]$ è per definizione la distanza $d_{T_\lambda}(x_\lambda,y_\lambda)$.
\end{defn}
\subsection{Complessi di bandiere}
I complessi di bandiere, o flag complexes, sono complessi simpliciali astratti che sono determinati dal loro $1$-scheletro, ovvero dal loro grafo sottostante. Ogni complesso simpliciale è quindi un sottocomplesso di un flag complex. 
\begin{defn}
	Sia $L$ un complesso simpliciale astratto e sia $V$ l'insieme dei suoi vertici. $L$ è un \textit{flag complex} se soddisfa la seguente condizione di \textit{assenza di triangoli}: \begin{center}
		Se $\{v_i,v_j\}\subseteq V$ è un simplesso di $L$, $\forall i,j\in \{1,\dots,n\}$, allora $\{v_1,\dots,v_n\}$ è un simplesso di $L$.
	\end{center}
\end{defn}
\begin{lem}{\label{flagLink}}
	Sia $K$ un flag complex. Allora il link in ogni vertice è un flag complex. 
	\begin{proof}
		Sia $L=\Lk(x,K)$. I vertici di $L$ sono gli $1$-simplessi di $K$ contenenti $x$. Siano $e_1,\dots,e_n$ vertici di $L$ tali che per ogni coppia di indici $i,j$, esista un $1$-simplesso $\sigma_{ij}$ di $L$ avente $e_i$ ed $e_j$ come estremi. Ciò significa che $\sigma_{ij}$ è un $2$-simplesso di $K$ contenente $e_i,e_j$. Posto $e_i=\{x,x_i\}$, si ha $\sigma_{ij}=\{x,x_i,x_j\}$ e dunque $\{x_i,x_j\}$ è un $1$-simplesso di $K$. \\
		Segue che nell'insieme $C=\{x,x_1,\dots,x_n\}$, ogni coppia di elementi costituisce un simplesso di $K$. Poiché $K$ è un flag complex, $C$ è un simplesso di $K$ contenente $x$ e ogni $e_i$. In conclusione, $C$ è un simplesso di $L$ contenente tutti i vertici $e_i$.
	\end{proof}
\end{lem}
\begin{defn}
	Un \textit{complesso sferico retto} è un complesso $M_1$-simpliciale le cui $1$-celle hanno tutte lunghezza $\pi/2$. 
\end{defn}
In un complesso sferico retto finito dimensionale è evidente che $\varepsilon(x)>0$ per ogni $x\in L$. 
\begin{lem}{\label{allrightLink}}
	Sia $L$ un complesso sferico retto. Allora il link in ogni punto è anch'esso un complesso sferico retto.
\end{lem}

\begin{oss}
	Sia $L$ un complesso simpliciale. Denotiamo con $\Phi_n(L)$ l'affermazione definita induttivamente nel modo seguente. \\
	$\Phi_0(L)$: esistono vertici $v_0,v_1,v_2\in L$ tali che $\{v_i,v_j\}$ è un simplesso $\forall i,j\in\{0,1,2\}$ ma $\{v_0,v_1,v_2\}$ non è un simplesso di $L$.\\
	$\Phi_{n+1}(L)$: esiste un vertcie $v_0$ tale che valga $\Phi_n(\Lk(v_0,L))$.\\
	Dalla definizione simpliciale di link, appare evidente che $L$ è un flag complex se e solo se $\Phi_n(L)$ è falsa per ogni $n> 0$.   
\end{oss}
\begin{thm}{\label{5.18}}
	Sia $L$ un complesso sferico retto finito dimensionale. $L$ è $\cat 1$ se e solo se è un flag complex.
	\begin{proof}
		Anzitutto osserviamo che se esistono vertici $v_0,v_1,v_2\in L$ congiunti a due a due da spigoli ma che non generano un simplesso, allora $[v_0,v_1]\cup[v_1,v_2]\cup[v_2,v_0]$ è una geodetica locale chiusa in $L$. Per dimostrarlo, notiamo che in $\Lk(v_1,L)$, che è un complesso sferico retto, i vertici corrispondenti agli spigoli $[v_0,v_1]$ e $[v_1,v_2]$ non sono congiunti da uno spigolo (il quale corrisponderebbe a un $2$-simplesso che li contiene). Quindi la distanza tra loro è almeno $\pi$ e quindi $[v_0,v_1]\cup[v_1,v_2]\cup[v_2,v_0]$ è una geodetica in un intorno di $v_1$\\ 
		Se $L$ è $\cat 1$, lo è anche il link in ogni vertice $v\in L$. Quindi sono $\cat 1$ anche il link di ogni vertice $v'\in \Lk(v,L)$, il link di ogni vertice $v''\in\Lk(v',\Lk(v,L))$, ecc. In particolare, nessuno di questi link può contenere un cerchio geodetico di lunghezza inferiore di $2\pi$. Ognuno di questi link è inoltre un complesso sferico retto. Quindi, deduciamo che l'affermazione $\Phi_n(L)$ è falsa $\forall n>0$, e quindi $L$ è un flag complex.\\
		Rimane da provare che se $L$ è un flag complex allora è uno spazio $\cat{1}$. 
		Procedendo per induzione sulla dimensione di $L$, possiamo assumere che $\Lk(v,L)$ sia $\cat 1$ per tutti i vertici $v$. Infatti, se $L$ è un flag complex di dimensione $n+1$, ogni link ha dimensione $n$ e dunque, per ipotesi induttiva, è $\cat{1}$ essendo un complesso sferico retto e flag, per i Lemmi \ref{flagLink} e \ref{allrightLink}. Mostriamo che ogni curva $\ell\subset L$ isometrica a un cerchio deve avere lunghezza almeno $2\pi$: segue dal Teorema \ref{5.4} che $L$ è $\cat{1}$. Sia $v$ un vertice di $L$ e sia $\ell'$ una componente connessa di $\ell\cap B_{\pi/2}(v)$ (qualora essa esista). Consideriamo lo sviluppo di $\ell'$ in $\mathbb S^2$: esso è una geodetica locale di lunghezza pari a quella di $\ell'$ e giace in un emisfero aperto di $\mathbb S^2$ cioè in $B_{\pi/2}(\bar v)$ (ove $\bar v$ è il punto della sfera che corrisponde a $v$). Poiché gli estremi di tale geodetica locale giacciono sul bordo dell'emisfero, essa deve avere lunghezza $\pi$. 
		\\
		Supponiamo dunque che la lunghezza di $\ell$ sia inferiore a $2\pi$. In tal caso $\ell$ non può contenere due archi disgiunti di lunghezza almeno $\pi$ e quindi non può intersecare due palle disgiunte di raggio $\pi/2$ centrate in vertici di $L$. Segue che se $\ell$ interseca due palle di raggio $\pi/2$ centrate in vertici, esiste un $1$-simplesso che congiunge tali vertici. Poiché $L$ è flag, l'insieme dei vertici $v\in L$ tali che $\ell$ interseca $B_{\pi/2}(v)$ genera un sotto-grafo completo di $L^{(1)}$, e dunque un simplesso $S$. Inoltre, $\ell$ sarebbe contenuta in $S$ e ciò è impossibile poiché nessuna sfera ha cerchi di lunghezza inferiore a $2\pi$. 
	\end{proof}
\end{thm}
\newpage
\chapter{Teorema di Gromov}
Il Teorema di Gromov, o criterio di Gromov per complessi cubici, consente di ridurre il problema di stabilire se un complesso ammetta una metrica di curvatura $\leq\ka$ ad un problema legato alla geometria dei link nei vertici del complesso.

\begin{defn}{\label{linkCn}}
	Un complesso $M_\ka$-poliedrico $K$ soddisfa la condizione dei link se $\Lk(v,K)$ è $\cat{1}$ per ogni vertice $v$ di $K$.
\end{defn}
\begin{prop}{\label{5.2}}
	Un complesso $M_\ka$-poliedrico $K$ con $\sh K$ finito ha curvatura $\leq \ka$ se e solo se soddisfa la condizione dei link.  
	\begin{proof}
		Sia $v$ un vertice di $K$. Per il Lemma \ref{7.39}., esiste $\varepsilon=\varepsilon(v)>0$ tale che $B_{\varepsilon}(v)$ sia convesso e isometrico alla palla di raggio $\varepsilon(v)$ centrata nel vertice di $C_\ka(\Lk(v,K))$. Per il Teorema  di Berestovskii (\ref{Ber}), $\Lk(v,K)$ è uno spazio $\cat 1$ se e solo se $C_\ka(\Lk(v,K))$ è $\catk$. Dunque, $K$ soddifa la link condition se e solo se ogni vertice ha un intorno che è $\catk$. 
		Dobbiamo dimostrare che $K$ ha curvatura $\leq \ka$ anche nei punti che non sono vertici. Questo discende dal fatto che ogni punto di $K$ ha un intorno isometrico all'intorno di un qualche punto sufficientemente vicino a un vertice del suo supporto (per il Lemma \ref{I.7.56}).
	\end{proof}
\end{prop}
\begin{exa}
	Sia $K$ un complesso $M_\ka$-poliedrico $2$-dimensionale. Il link di un vertice $v\in K$ è un grafo metrico (non necessariamente connesso) i cui vertici corrispondono alle $1$-celle incidenti in $v$ e i cui archi corrispondono agli angoli delle $2$-celle incidenti in $v$. La lughezza di un arco è uguale all'ampiezza dell'angolo interno in $v$ della $2$-cella corrispondente.
\end{exa}
Possiamo quindi dimostrare il Teorema di Gromov.
\begin{thm}[Criterio di Gromov]{\label{Gromov}}
	Un complesso cubico finito dimensionale ha curvatura non-positiva se e solo se il link di ognuno dei suoi vertici è un complesso di bandiere.
	\begin{proof}
		Per la Proposizione \ref{5.2}, un complesso cubico (che in particolare è un complesso euclideo, o $M_0$-complesso) ha curvatura non-positiva se e solo se il link in ogni vertice è $\cat 1$. Un tale link è un complesso sferico retto, quindi possiamo applicare il Teorema \ref{5.18}, concludendo che un complesso cubico ha curvatura $\leq 0$ se e solo se ogni link è $\cat 1$.
	\end{proof}
\end{thm}
Segue il seguente risultato fondamentale.
\begin{cor}
	Un complesso cubico finito dimensionale e connesso è uno spazio $\cat 0$ se e solo se è semplicemente connesso e i link di ognuno dei suoi vertici è un complesso di bandiere. 
	\begin{proof}
		Segue dal Criterio di Gromov (\ref{Gromov}) e dal Teorema \ref{5.4}.
	\end{proof}
\end{cor}
\newpage

\chapter{Right-Angled Artin Groups Orientati}
Nell'ultima sezione introduciamo la nozione di RAAG (Right-Angled Artin Group) orientato che generalizza quella di RAAG definito su un grafo. Dimostreremo che ad un tale gruppo è possibile associare in modo naturale uno spazio $\cat 0$ su cui esso agisce geometricamente. Investigheremo alcune proprietà elementari di questa classe di gruppi stabilendo alcune relazioni con i gruppi di Artin e di Coxeter. 
\section{Grafi e gruppi} 
In letteratura esistono numerose definizioni di grafo, ognuna delle quali risulta conveniente in un particolare ambito. Per esempio, in \cite{S}, un grafo è una coppia di insiemi $(V,E)$ munita di due applicazioni $\bar{}:E\to E$ e $(o,t):E\to V^2$ tali che $e\mapsto\bar e$ è una involuzione non banale e $(o,t)(\bar e)=(t,o)(e)$. Nel seguito, useremo essenzialmente due definizioni di grafo che si adattano meglio al nostro contesto. \\
Un \textit{grafo na\"{i}ve} è una coppia $\Gamma_n=(\bar V,\bar E)$ ove $\bar V$ è un insieme finito\footnote{La definizione di grafo non richiede che l'insieme dei vertici sia finito; tuttavia, per semplicità, limiteremo il nostro studio a questo caso.} non vuoto ed $\bar E\subseteq \mathcal P_2(\bar V)$ è un insieme di sottoinsiemi di cardinalità $2$ di $\bar V$.\\
Un \textit{grafo orientato} è una coppia $\Gamma=(V,E)$ ove $V$ è un insieme finito non vuoto di vertici ed $E\subseteq V\times V\setminus \Delta(V)$ è un insieme di coppie ordinate di vertici distinti di $V$. Per ogni arco $e=(v,w)\in E$, possiamo considerare le projezioni sull'insieme dei vertici. In particolare, poniamo $o(e)=v,\ t(e)=w$.\\
Ad un grafo orientato $\Gamma=(V,E)$ possiamo associare in modo canonico un grafo na\"ive $\ddot\Gamma=(\ddot V,\ddot E)$ ponendo: 
\begin{enumerate}
	\item $\ddot V=V$,
	\item $\ddot E\subseteq \mathcal P_2(\ddot V)$: $\{v,w\}\in \ddot E\ \Leftrightarrow\ \{(v,w),(w,v)\}\cap E\neq \emptyset$. 
\end{enumerate}
Chiaramente, le funzioni $o$ e $t$ non si inducono sul grafo na\"ive associato.\\

Esiste una classe ristretta di grafi orientati su cui ci concentreremo nel seguito. Diremo che $\Gamma=(V,E)$ è un \textit{grafo speciale} se l'insieme dei vertici $V$ (risp. l'insieme degli archi $E$) ammette una decomposizione $V=V_0\cup V_s$  (risp. $E=E_0\cup E_s$) in vertici (risp. archi) normali e speciali e un vertice è speciale solo se è il punto finale di un arco speciale. Richiediamo inoltre che ogni clique di $\Gamma$ contenga al più un vertice speciale.\\
Un grafo speciale può essere pensato come un grafo orientato avente un solo arco per ogni arco speciale (orientato nello stesso modo) e una coppia di archi mutuamente inversi (ovvero un arco geometrico, o non orientato) per ogni arco normale. Possiamo inoltre supporre che un arco sia speciale solo se il suo vertice finale è speciale. Un vertice del grafo è detto \textit{attrattore} se ogni arco che lo contiene è speciale (e quindi il vertice è punto finale di ogni tale arco).
\begin{defn}
	Sia $\Gamma=(V,E)=(V_0\cup V_s,E_0\cup E_s)$ un grafo speciale. Definiamo il \textit{RAAG orientato}\footnote{Nel seguito, chiameremo RAAG sia i right-angled Artin groups sia i RAAG orientati. Quando il grafo non è orientato, chiameremo il gruppo ad esso associato RAAG classico.} associato a $\Gamma$ tramite la seguente presentazione:
	\[G_\Gamma:=\langle V\ \vert \ R(e)\ :e\in E\rangle \]
	ove $R(e)$ sono le relazioni definite nel seguente modo: \begin{enumerate}
		\item Se $e\in E_0$, $R(e)=[t(e),o(e)]=t(e)o(e)t(e)^{-1}o(e)^{-1}$
		\item Se $e\in E_s$, $R(e)=[\!\![t(e),o(e)]:=t(e)o(e)t(e)^{-1}o(e)$.
	\end{enumerate}
I generatori corrispondenti ai vertici di $\Gamma$ e i loro inversi sono i \textit{generatori canonici} di $G_\Gamma$.

\end{defn}
Osserviamo che la relazione associata a un arco speciale è esattamente quella che si impone sui generatori naturali del gruppo fondamentale della bottiglia di Klein; per questo motivo, diamo ad essa il nome di \textit{relazione di Klein}. Il gruppo fondamentale della bottiglia di Klein è il più semplice gruppo di Artin orientato che non è un gruppo di Artin classico.\\
Dimostriamo ora una identità che vale per gli elementi di un gruppo che soddisfano la relazione di Klein. Questa verrà largamente usata nel seguito.
\begin{lem}{\label{idKl}}
	Sia $G$ un gruppo e siano $x,y\in G$ tali che $[\!\![x,y]=xyx^{-1}y=1$. Allora vale $\forall a,b\in \mathbb Z$:
	\[x^ay^b=y^{(-1)^ab}x^a.\]
	\begin{proof}
		Se $a=0$, non c'è nulla da dimostrare. Sia quindi $a>0$. Poiché $xyx^{-1}=y^{-1}$ ricaviamo per induzione: \begin{align*}
		x^ay^b&=x^{a-1}(xy^bx^{-1})x=\\
	    &=x^{a-1}(xyx^{-1})^bx=\\
		&=x^{a-1}y^{-b}x=\dots=\\
		&=y^{(-1)^ab}x^a
		\end{align*}
		Se $a<0$ si ragiona in modo analogo.
	\end{proof}
\end{lem}

\begin{oss}{\label{isograph}}
	\`E ben noto (per esempio \cite{Dr}) che dati due gruppi right-angled di Artin (o di Coxeter) definiti da grafi $\Gamma_1$ e $\Gamma_2$, se $\Gamma_1\neq\Gamma_2$, anche i gruppi sono diversi. Nel caso orientato questo non vale, come mostra l'esempio seguente.\\
	Consideriamo i due grafi con $3$ vertici $\{v_1,v_2,v_3\}$ e archi rispettivamente $\{(v_1,v_2),(v_1,v_3),(v_3,v_1)\}$ e $\{(v_1,v_2),(v_1,v_3)\}$.  Graficamente,
	
	$\Gamma_1=$ \begin{tikzcd}
	v_2                              &     \\
	v_1 \arrow[u] \arrow[r, no head] & v_3
	\end{tikzcd}\hspace{3cm}
	$\Gamma_2=$
	\begin{tikzcd}
	v_2                     &     \\
	v_1 \arrow[u] \arrow[r] & v_3
	\end{tikzcd}\\
	Chiaramente, i due grafi sono speciali e distinti. Siano $G_i=G_{\Gamma_i}$, $i=1,2$. I due gruppi sono \begin{align*}
	G_1&=\langle a,b,c\ \vert\ [\!\![b,a],[a,c]\rangle\\
	G_2&=\langle x,y,z\ \vert\ [\!\![y,x],[\!\![z,x]\rangle.
	\end{align*}
	L'applicazione $a\mapsto x, b\mapsto y, c\mapsto zy^{-1}$ è ben definita poiché
	\begin{align*}
	[zy^{-1},x]&=(zy^{-1})x(zy^{-1})^{-1}x^{-1}=\\
	&=zy^{-1}xyz^{-1}x=\ \ \ \ \ \ \ \ \ \   (x^{-1}=y^{-1}xy)\\
	&=zx^{-1}z^{-1}x^{-1}=1\ \ \ \ \ \ \ \ \ \  (zxz^{-1}x=1)
	\end{align*}
	Inoltre, essa è chiaramente un isomorfismo perché ammette inversa $x\mapsto a, y\mapsto b, z\mapsto cb$. Quindi si ha $G_1\simeq G_2$, anche se $\Gamma_1\not\simeq \Gamma_2$.
	\section{Alcune proprietà dei RAAG}
	Come vedremo più avanti, sotto opportune ipotesi sul grafo speciale, il RAAG orientato ad esso associato è un prodotto amalgamato dei gruppi associati alle clique massimali del grafo. Studiamo quindi tali gruppi \textit{locali} concentrandoci sul rapporto che sussiste tra il gruppo e il RAAG classico associato alla clique.\\
	Sia $\Lambda$ una clique speciale con un vertice speciale. Supponiamo che $\Lambda$ sia un cono con un attrattore e base non-orientata; in modo equivalente $\Lambda$ è un grafo completo e ha lo stesso numero di archi speciali e vertici normali. Esplicitamente, $V=\{v_s,v_1,\dots,v_n\}$, con $V_s=\{v_s\}$ e $E_s=\{(v_1,v_s),\dots,(v_n,v_s)\}$. Allora \[G_\Lambda=\langle a,b_1,\dots,b_n\ \vert\ [\!\![a,b_i],[b_i,b_j]:\ i,j=1,\dots,n\rangle.\]
	Vale il seguente:
	\begin{lem}{\label{semidir}}
		Nelle ipotesi precedenti, \[G_\Lambda\simeq \mathbb Z^n\rtimes _\tau\mathbb Z,\]
		ove $\tau:\mathbb Z\to \operatorname{Aut}(\mathbb Z^n)$ è definito da $\tau(z)(w)=(-1)^zw$.\\
		Quindi, se $\Lambda$ è un grafo completo speciale con $n$ archi speciali e $m$ archi normali contenenti il vertice speciale, \[G_\Lambda\simeq \mathbb Z^m\times(\mathbb Z^n\rtimes _\tau\mathbb Z)\]
		\begin{proof}
			Sia $f:G_\Lambda\to \mathbb Z^n\rtimes _\tau\mathbb Z$ l'unico morfismo di gruppi tale che $f(a)=(0,1)$ e $f(b_i)=(e_i,0)$, ove $\{e_1,\dots,e_n\}$ è la base standard di $\mathbb Z^n$. La funzione $f$ esiste perché valgono le seguenti identità $\forall i,j=1,\dots,n$: \begin{align*}
			[\!\![f(a),f(b_i)]=&(0,1)(e_i,0)(0,1)^{-1}(e_i,0)=\\
			=& (0,1)(e_i,0)(0,-1)(e_i,0)=\\
			=&(-e_i,1)(-e_i,-1)=(-e_i+(-1)^{1}(-e_i),1-1)=\\
			=&(0,0)
			\end{align*} 
			e
			\begin{align*}
			[f(b_i),f(b_j)]=&(e_i,0)(e_j,0)(e_i,0)^{-1}(e_j,0)^{-1}=\\
			=& (e_i,0)(e_j,0)(-e_i,0)(-e_j,0)=\\
			=& (0,0).		\end{align*}
			Il morfismo $f$ è un isomorfismo perché ammette inversa $g:\mathbb Z^n\rtimes _\tau\mathbb Z\to G_\Lambda$ definita da \[(v,z)\longmapsto b_1^{v_1}b_2^{v_2}\dots b_n^{v_n}a^k.\]
			$g$ è un morfismo di gruppi perché $\forall v,v'\in\mathbb Z^n,\ \forall k,k'\in \mathbb Z$ si ha, grazie al Lemma \ref{idKl} e al fatto che i $b_i$ commutano: \begin{align*}
			g(v,k)g(v',k')&=b_1^{v_1}b_2^{v_2}\dots b_n^{v_n}a^kb_1^{v'_1}b_2^{v'_2}\dots b_n^{v'_n}a^{k'}=\\
			&=b_1^{v_1}\dots b_n^{v_n}b_1^{(-1)^k v'_1}a^kb_2^{v'_2}\dots b_n^{v'_n}a^{k'}=\\
			&=\dots=\\
			&=b_1^{v_1}\dots b_n^{v_n}b_1^{(-1)^k v'_1}\dots b_n^{(-1)^kv_n}a^ka^{k'}=\\
			&=b_1^{v_1+(-1)^kv'_1}\dots b_n^{v_n+(-1)^kv'_n}a^{k+k'}=\\
			&=g(v+(-1)^kv,k+k')=g((v,k)(v',k')).
			\end{align*}
			Infine, una clique con $n$ archi speciali e $m$ archi normali è il join (non orientato) di una clique speciale $\Lambda_s$ con $n$ archi speciali con una clique non-orientata $\Lambda_0$ su $m$ vertici. Dunque il gruppo ad essa associato è il prodotto diretto \[G_\Lambda=G_{\Lambda_s}\times G_{\Lambda_0}.\]
		\end{proof}
	\end{lem}  
	\begin{cor}
		Sia $\Lambda$ una clique speciale con lo stesso numero di archi speciali e vertici normali. Allora il centro di $G_\Lambda$ è \[Z(G_\Lambda)=\langle a^2\rangle\simeq \{0\}\times 2\mathbb Z, \] ove $a$ è il generatore di $G_\Lambda$ corrispondente al vertice speciale di $\Lambda$.  
		\begin{proof}
			Sia $G=G_\Lambda=\mathbb Z^n\rtimes \mathbb Z$ come nel Lemma \ref{semidir}. Calcoliamo il centro di $G$. Sia $(v,k)\in Z(G)$: allora $\forall (v',k')\in G$, \[(v,k)(v',k')=(v+(-1)^kv',k+k')=(v'+(-1)^{k'}v,k+k')=(v',k')(v,k).\] L'identità vale se e solo se $v=0$ e $k\in 2\mathbb Z$.
		\end{proof}
	\end{cor}
	\begin{prop}
		Sia $\Lambda$ una clique speciale. Sia $A_\Lambda$ il RAAG classico associato a $\Lambda$. Allora $A_\Lambda\simeq \mathbb Z^{\vert V\vert}$ si immerge come sottogruppo normale di indice finito in $G_\Lambda$.\begin{proof} Se $\Lambda$ non ha vertici speciali, $G_\Lambda=A_\Lambda$ è un gruppo abeliano libero.\\
			Supponiamo quindi che esista un vertice speciale in $\Lambda$. Siano $n$ il numero di archi speciali di $\Lambda$ e $m$ il numero di archi normali contenenti il vertice speciale. Identifichiamo $G_\Lambda$ con il gruppo $\mathbb Z^m\times(\mathbb Z^n\rtimes _\tau\mathbb Z)$. \`E chiaro che $A_\Lambda=\mathbb Z^{m+n+1}=\mathbb Z^m\times \mathbb Z^n\times \mathbb Z$. Consideriamo l'applicazione $f:A_\Lambda\to G_\Lambda$ che associa a un generatore canonico $a$ di $A_\Lambda$ il suo quadrato $a^2$ in $G_\Lambda$. Essa è ben definita. Inoltre, se $f(w,v,k)=0$, vale $(2w,2v,2k)=0$, cioè $(w,v,k)=0$; questo prova che $f$ è injettiva. \\
			Sia poi $\pi:G_\Lambda\to W_\Lambda$ la projezione sul gruppo di Coxeter associato a $\Lambda$. Si ha che $f(A_\Lambda)$ è il nucleo di $\pi$, ovvero \[1\to A_\Lambda\to G_\Lambda\to W_\Lambda\to 1\] è una successione esatta.\\
			Poiché $W_\Lambda$ è finito, $A_\Lambda$ ha indice finito in $G_\Lambda$.
		\end{proof}
	\end{prop}

	
\end{oss}
\section{Bottiglie di Klein n-dimensionali e coomologia}
Abbiamo già notato che il gruppo di Artin orientato associato ad un segmento con un attrattore è il gruppo fondamentale della bottiglia di Klein. In questa sezione generalizziamo questa situazione al caso di un grafo completo con un vertice speciale, costruendo una bottiglia di Klein di dimensione arbitraria. Questo ci consentirà di costruire uno spazio classificante per i RAAG che generalizza il complesso di Salvetti. Il complesso così costruito ci consentirà di determinare la coomologia del gruppo, almeno per una classe di grafi soggiacenti.
\begin{defn}
	Sia $m$ un numero intero positivo. Poniamo \[K_m:=(S^1)^m/\sim\]
ove $\sim $ è la relazione di equivalenza generata da $(z_1,\dots,z_{m-1},z_m)\sim (\bar z_1,\dots,\bar z_{m-1},-z_m)$. $K_m$ munito della topologia quoziente è detto \textit{bottiglia di Klein $m$-dimensionale}: per $m=2$, $K_2$ è la bottiglia di Klein .
\end{defn} 
Nell'articolo \cite{D} viene dimostrato che il gruppo fondamentale\footnote{Data la connessione per archi degli spazi considerati, il punto base non compare esplicitamente nella notazione di gruppo fondamentale: ogni gruppo fondamentale considerato sarà dunque definito a meno di isomorfismi.} di $K_m$ è  \[\pi_1(K_m)\simeq \langle a_1,\dots,a_{m-1},b\ \vert\ ba_ib^{-1}a_i,\ [a_i,a_j]:\ i,j=1,\dots,m-1\rangle.\]
Dunque il gruppo fondamentale della bottiglia di Klein $(n+1)$-dimensionale è il RAAG orientato sul grafo completo $\Lambda$ con $n$ vertici normali e un attrattore. 
Inoltre, si vede subito che il rivestimento universale di $K_{n+1}$ è $\mathbb R^{n+1}$, cioè $K_{n+1}$ è un $K(G_\Lambda,1)$.
\begin{thm}
	Sia $G$ il RAAG orientato associato ad un grafo completo con $n$ vertici normali e un attrattore. Allora la bottiglia di Klein $(n+1)$-dimensionale è un $K(G,1)$. 
\end{thm}

Definiamo ora \[K_{m,l}=(S^1)^m/\sim'\]
ove $\sim' $ è la relazione di equivalenza generata da $(z_1,\dots,z_m)\sim (w_1,\dots,w_m)$ con $w_i=\bar z_i$ se $i\leq l$ e $w_i=-z_i$ se $i>l$. Vale $$K_{m,l}\approx K_{l+1}\times (S^1)^{m-l-1}.$$ Infatti, l'applicazione $h:K_{m,l}\to K_{l+1}\times (S^1)^{m-l-1}$ definita da \[[(z_1,\dots,z_m)]\mapsto ([(z_1,\dots,z_{l+1})],(z_{l+1}^{-1}z_{l+2},\dots,z_{l+1}^{-1}z_{m}))\] è continua,
 e ammette inversa continua \[\left([(z_1,\dots,z_{l+1})],(z_{l+2},\dots,z_m)\right)\mapsto [z_1,\dots,z_{l+1},z_{l+1}z_{l+2},\dots,z_{l+1}z_{m}]\]
 
Segue che il suo gruppo fondamentale è il prodotto diretto dei gruppi fondamentali dei suoi fattori, i.e. $$\pi_1(K_{m.l})= \pi_1(K_{l+1})\times \pi_1((S^1)^{m-l-1})$$ ed è quindi (isomorfo a) il RAAG orientato definito dal grafo completo $\Lambda$ con $l$ archi speciali e $m$ vertici.\\
Osserviamo che $K_{m,l}$ ha una naturale decomposizione in celle che lo rende un CW-complesso. \\

Ancora in \cite{D} si dimostra la seguente formula:
\begin{thm}
	Sia $K_m$ la bottiglia di Klein $m$-dimensionale. Allora il suo anello di coomologia a coefficienti nel campo con $2$ elementi è (isomorfo a):
	$$H^\bullet(K_m,\mathbb F_2)=\frac{\mathbb F_2[R,V_1,\dots,V_m]}{(R^2,V_i^2+RV_i)},$$ ove il grado dei generatori $R$ e $V_i$ è pari a $1$.
\end{thm}
Dalla formula di K\"unneth segue che \begin{align*}
H^\bullet (K_{m,l},\mathbb F_2)=& \frac{\mathbb F_2[R,V_1,\dots,V_l]}{(R^2,V_i^2+RV_i)}\otimes \Lambda_{\mathbb F_2}(W_1,\dots,W_{m-l-1})
\end{align*}
cioè \begin{align}\begin{split}H^\bullet (K_{m,l},\mathbb F_2)=& \frac{\mathbb F_2[R,V_1,\dots,V_l]}{(R^2,V_i^2+RV_i)}\otimes \frac{\mathbb F_2[W_1,\dots,W_{m-l-1}]}{(W_i^2)}=\\
=&\frac{\mathbb F_2[R,V_1,\dots,V_l,W_1,\dots,W_{m-l-1}]}{(R^2,V_i^2+RV_i,W_i^2)}\end{split}\end{align}
ove $\vert R\vert=\vert V_i\vert =\vert W_j\vert=1$.\\
Se $\Lambda=(V,E)$ è un grafo speciale completo con $\vert E_s\vert =r$ e $\vert V\vert =n+1$, essendo $K_{n,r}$ un $K(G_\Lambda,1)$ vale:
\[H^\bullet (G_\Lambda,\mathbb F_2)=\frac{\mathbb F_2[R,V_1,\dots,V_l]}{(R^2,V_i^2+RV_i)}\otimes \frac{\mathbb F_2[W_1,\dots,W_{m-l-1}]}{(W_i^2)}.\]

\subsection{Algebra associata a un grafo speciale}
Sia $\Gamma$ un grafo speciale. Definiamo una $\mathbb F_2$-algebra $\mathbb F_2(\Gamma)$ associata a $\Gamma$. Essa è l'anello di polinomi nelle variabili $v$ in corrispondenza biunivoca con i vertici di $\Gamma$ soggetto alle seguenti relazioni:
\begin{enumerate}
	\item se $v\in V$ non è il punto inziale di un arco speciale (per esempio se $v\in V_s$), $v^2=0$;
	\item se $w\in V_0$ è estremo di un arco speciale $(w,v)\in E_s$, $vw+w^2=0$;
	\item se $\{v,w\}\notin \ddot E$, $vw=0$.
	\item se $\{v,w\}\notin \ddot E$ ed esiste $v'\in V$ tale che $(v',v)\in E_s$, $(v',w)\in E$, si impone $vv'+wv'=0$.
\end{enumerate}
$\mathbb F_2(\Gamma)$ è un'algebra quadratica e coincide con l'algebra associata a $\ddot \Gamma$ (\cite{Dr}) nel caso in cui $\Gamma$ non abbia vertici speciali. 
\begin{exa}
	Sia $\Gamma$ non orientato, i.e. $V=V_0$, $\ddot E:=E=E_0$. L'anello di coomologia di $G_\Gamma$ è il quoziente della $\mathbb F_2$-algebra esterna generata da $V$ rispetto all'ideale generato dai monomi $vw$ per ogni $\{v,w\}\notin E$. Dato che stiamo lavorando in caratteristica $2$, essa coincide con $\mathbb F_2(\Gamma)$.
\end{exa}
\newcommand{\f}{{\mathbb{F}_2}}
\begin{exa}
	Sia $\Gamma$ un grafo speciale completo con $n+1$ vertici. Anche in questo caso, la coomologia di $G_\Gamma$ coincide con $\mathbb F_2(\Gamma)$. Infatti, se $V_s=\{s\}$ e $\vert E_s\vert=r$, l'algebra associata a $\Gamma$ è 
	\[\frac{\f[s,v_1,\dots,v_r,w_1,\dots,w_{n-r}]}{(sv_i+v_i^2,s^2,w_j^2)}\]
\end{exa}
Dunque, $\Gamma\to \f(\Gamma)$ è un buon candidato per la coomologia di $G_\Gamma$.
\begin{oss}
	Le algebre associate ai due grafi dell'Osservazione \ref{isograph} con gruppi isomorfi sono isomorfe. Siano $\Gamma_1$, $\Gamma_2$ come nell'Osservazione. Sia $H_i=H^\bullet(G_{\Gamma_i},\mathbb F_2)$ ($i=1,2$). Allora \begin{align*}
	H_1=\frac{\mathbb F_2[x,y,z]}{(x^2,z^2,xz,xy+y^2,xy+zy)}\\
	 H_2=\frac{\mathbb F_2[x,y,z]}{(x^2,z^2,xz,xy+y^2,zy+y^2)}
	\end{align*} Poiché in $H_2$ valgono $y^2=yz=xy$, vale anche $xy+yz=0$. Analogamente, dato che in $H_1$ valgono $xy=yz=y^2$, vale anche $yz+y^2=2y^2=0$. I due anelli sono quindi identici.
\end{oss}
\begin{conj}
	L'anello di coomologia a coefficienti in $\mathbb F_2$ di $G_\Gamma$ è isomorfo a $\mathbb F_2(\Gamma)$.
\end{conj}


\section{Complesso di Klein-Salvetti}
Abbiamo visto che le bottiglie di Klein $n$-dimensionali sono degli spazi classificanti per i RAAG associati a grafi speciali completi. In questa sottosezione generalizziamo le considerazioni fatte precedentemente e costruiamo un complesso di celle che è un $K(G_\Gamma,1)$ per un RAAG orientato arbitrario $G_\Gamma$, in modo simile al complesso di Salvetti associato a un RAAG classico. Come nel caso del complesso di Salvetti, eseguiamo due costruzioni differenti che conducono al medesimo complesso.\\

Nel seguito, sia $\Gamma=(V,E)$ un grafo speciale fissato e $G=G_\Gamma$ il RAAG orientato ad esso associato.\\
Richiamiamo la definizione di \textit{complesso ipercubico associato al gruppo} $G$ con sistema di generatori $V$ (Cfr. \cite{C}). Consideriamo il grafo di Cayley $\cay(G,V)$ di $G$ con sistema di generatori $V$. Esplicitamente, esso ha $G$ come insieme di vertici e le coppie $(v,gv^{\pm})$, con $v\in V,\ g\in G$, come archi orientati. Il complesso ipercubico $\mathcal C$ è determinato da $\cay(G,V)$ nel seguente modo: se $2^{n-1}n$ archi sono l'$1$-scheletro di un $n$-cubo, attacchiamo un tale $n$-cubo a $\cay(G,V)$ identificandone l'$1$-scheletro con gli archi considerati. L'$1$-scheletro di $\mathcal C$ è dunque il grafo di Cayley $\mbox{Cay}(G,V)$ di $G$ e il $2$-scheletro è un complesso di Cayley per $G$ (Cfr.\cite{hat}).
\\
\textbf{$1^\circ$ costruzione.} Sia $\mathcal C$ il \textit{complesso ipercubico} associato a $(G,V)$.  L'azione di $G$ sul suo graof di Cayley induce un'azione cellulare su e risulta quindi definito il quoziente $B=G\backslash \mathcal C$. Poiché $G$ agisce transitivamente su se stesso, $B$ ha un solo vertice $x_0$. Inoltre, per ogni vertice $v$ del grafo $\Gamma$, l'$1$-scheletro di $B$ presenta un cappio $e_v$ in $x_0$. Il $2$-scheletro di $B$ è ottenuto nel seguente modo: ad ogni arco $e\in E$, si attacca una $2$-cella al bouquet di cerchi $B^{(1)}$ identificando il suo bordo come prescritto dalla relazione $R(e)$.\\
\textbf{$2^\circ$ costruzione.} L'idea della costruzione è di attaccare delle bottiglie di Klein in un modo prescritto dal grafo $\Gamma$, così come il complesso di Salvetti è ottenuto attacando tori di dimensione crescente con opportune identificazioni.\\ Definiamo $B_\Gamma$ in modo induttivo. $B_\Gamma^{(0)}$ consiste di un solo punto $x_0$. $B_\Gamma^{(1)}$ ha un cappio in $x_0$ per ogni vertice $v$ del grafo $\Gamma$. Costruiamo lo scheletro di dimensione $n+1$ di $B_\Gamma$ come segue: per ogni $(n+1)$-clique $\Lambda$ di $\Gamma$, consideriamo il complesso $K_\Lambda=K_{n+1,r}=K_{r+1}\times (S^1)^{n-r}$ ove $r$ è il numero di archi speciali di $\Lambda$, e attribuiamo un'etichetta ad ogni sua coordinata in modo che le prime $r+1$ coordinate corrispondano ai vertici che giacciono su un arco speciale. Attacchiamo poi $K_\Gamma$ a $B_\Gamma^{(n)}$ identificando il complesso associato a una $n$-clique di $\Lambda$ con il sottocomplesso di $K_\Gamma$ generato dai vertici corrispondenti alla sottoclique. In questo modo $B_\Gamma$ è un complesso di celle ottenuto attaccando tori e bottiglie di Klein.
\begin{thm}
	$B_\Gamma$ è il quoziente $G\backslash \mathcal C$ ed è un spazio classificante di $G$. Inoltre, $\mathcal C$ è uno spazio $\cat 0$.
	\begin{proof}
		Mostriamo anzitutto che $\mathcal C$ è un complesso cubico $\cat 0$. 
Poiché $\mathcal C$ ha come $2$-scheletro un complesso di Cayley per $G$ e questo determina il suo gruppo fondamentale (per il Teorema di Seifert-Van Kampen), $\mathcal C$ è semplicemente connesso. Per il Teorema \ref{5.4}, resta da dimostrare che $\mathcal C$ soddisfa la condizione dei link.\\
In forza del Criterio di Gromov \ref{Gromov}, basta mostrare che i link nei vertici di $\mathcal C$ sono complessi di bandiere. Siano quindi $e_1,\dots,e_n$ archi in $\cay(G,V)$ uscenti dal vertice corrispondente a $1\in G$, i.e. $e_i\in \Lk(1,\mathcal C)$. Ogni arco è della forma $e_i=(1,x_i)$ per qualche $x_i\in V^{\pm}= V\coprod V^{-1}=\{v^{\pm}\vert\ v\in V\}$. Dobbiamo dimostrare che, se $\forall i<j$, $e_i,e_j$ sono contenuti in una $2$-cella di $\mathcal C$ (cioè in un quadrato), $e_1,\dots,e_n$ sono contenuti in un $n$-cubo. Poiché $e_i$ ed $e_j$ sono contenuti in un quadrato, $\{1,x_i,g_{ij},x_j\}$ è l'insieme di vertici di un quadrato in $\cay(G,V)$. 
 L'elemento $g_{ij}\in G$ è adiacente in $\cay(G,V)$ sia a $x_i$ sia a $x_j$, dunque esistono generatori $x_i',x_j'\in V^{\pm}$ tali che $g_{ij}=x_ix_j'$ e $g_{ij}x_i'=x_j$. Si ricava $x_ix_j'x_i'=x_j$. Chiaramente, non può essere $x_j'=x_i^{-1}$, altrimenti, $g_{ij}=x_ix_j'=x_ix_i^{-1}=1$ e $\{1,x_i,g_{ij},x_j\}$ non definisce un quadrato.
Se $i<j$ si ha una delle seguenti:
\begin{align}
\label{prima}x_i'=&x_i^{-1}\mbox{, } x_j'=x_j,\ \mbox{oppure}\\
\label{seconda} x_i'=&x_i,\ x_j'=x_j.
\end{align} Se vale (\ref{prima}), \[g_{ij}=x_ix_j\mbox{ e }[x_i,x_j]=1.\] Considerando la projezione sul RACG associato a $\Gamma$, si vede che la commutazione dei generatori canonici implica che $(v_i,v_j)\in E_0$, ove $v_l$ è il vertice corrispondente al generatore canonico $x_l$.\\
Se invece vale (\ref{seconda}), \[g_{ij}=x_ix_j\mbox{ e }\kl x_j,x_i]=1.\] Considerando ancora la projezione, si vede che necessariamente $(v_i,v_j)\in E_s$.\\ 
Segue che nel grafo $\Gamma$ che definisce $G$, i vertici corrispondenti ai generatori $x_1,\dots,x_n$ sono a due a due adiacenti e quindi \[\{1,\ x_{i_1}\dots x_{i_k}\vert\ {i_1<\dots<i_k},\ 1\leq k\leq n\}\] è l'insieme di vertici di un $n$-cubo. Concludiamo che gli archi $e_1,\dots,e_n$ sono contenuti in un $n$-cubo, ovvero, come vertici di $\Lk(1,\mathcal C)$, generano un $(n-1)$-simplesso.\\

Mostriamo ora che $B_\Gamma=G\backslash \mathcal C$. Anzitutto, se $\Gamma$ è un grafo speciale completo (con o senza vertice speciale), $G_\Gamma=\mathbb Z^m\times (\mathbb Z^n\rtimes_\tau \mathbb Z)$ ove $n$ è il numero di archi speciali, $m+n+1=\vert V\vert$ è il numero di vertici e $\tau$ è il morfismo banale se $\Gamma$ è priva di vertice speciale; se $\Gamma$ ha un vertice speciale, $\tau:\mathbb Z\to \mbox{Aut}\mathbb Z^n$ è definito da $\tau(k)(v)=k\cdot v=(-1)^kv$ come nel Lemma \ref{semidir}. Inoltre, $\mathcal C$ omeomorfo a $\mathbb R^{m+n+1}$. L'azione di $G$ su $\mathcal C$ è definita da \[(w,v,k)(x,y,z)=(x+w,(-1)^ky+v,z+k),\] per $(w,v,k)\in  \mathbb Z^m\times (\mathbb Z^n\rtimes_\tau \mathbb Z)$ e $(x,y,z)\in \mathbb R^m\times\mathbb R^n\times\mathbb R$. Segue che $G\backslash\mathcal C$ è omeomorfo a $K_{m+n+1,n}=K_\Gamma=B_\Gamma$. 
	\end{proof}
\end{thm}
\begin{cor}
	Sia $\Gamma$ un grafo speciale $\operatorname{(finito)}$. Allora $G_\Gamma$ è privo di torsione.
	\begin{proof}
		Poiché $G_\Gamma$ ha uno spazio classificante di dimensione finita, è privo di torsione (\cite{hat}).
	\end{proof}
\end{cor}
\section{Grafo di gruppi associato a un RAAG}
Concludiamo dimostrando che, sotto opportune ipotesi sul grafo soggiacente, il RAAG può essere decomposto come gruppo fondamentale di un albero di gruppi (nel senso di \cite{S}).
\begin{defn}
 Un grafo orientato $\Gamma$ si dice \textit{cordale} se ogni cammino in $\ddot \Gamma$ di lunghezza almeno $4$ ha una corda, i.e. per ogni $n$-pla $(v_1,\dots, v_n)$ di vertici tali che $\{v_i,v_{i+1}\}\in \ddot E$, $i=1,\dots,n-1$, se $n>3$, esistono indici $i<j\neq i+1$ tali che $\{v_i,v_j\}\in \ddot E$.   
\end{defn}
\newcommand{\cli}{\mathcal{K}_\Gamma}
Dato un grafo $\Gamma$, sia $\mathcal K_\Gamma$ l'insieme delle clique massimali di $\Gamma$. Un \textit{albero di clique} di $\Gamma$ è un albero (nel senso di \cite{S}) con insieme di vertici $\cli$.\\
Diciamo che un albero di clique $T$ di $\Gamma$ soddisfa la \textit{proprietà dell'intersezione} se per ogni clique massimali $\Lambda,\Lambda'\in\cli$, l'intersezione $\Lambda\cap\Lambda'$ è contenuta in ogni clique massimale del cammino minimo in $T$ da $\Lambda $ a $\Lambda'$.\\
Vale il seguente: 
\begin{thm}\cite[Thm. 3.2]{BP}
	Un grafo $\Gamma$ è cordale se e solo se esiste un albero di clique $T$ per $\Gamma$ che soddisfa la proprietà dell'intersezione. 
\end{thm}
Siamo pronti per definire un grafo di gruppi associato a un grafo speciale cordale $\Gamma$. Sia $T$ un albero di clique di $\Gamma$ con la proprietà di intersezione. Ad ogni clique massimale $\Lambda\in \cli$ associamo il gruppo \[G_\Lambda=\langle V\Lambda\ \vert \ R(e):\ e\in E\Lambda\rangle;\] ad ogni arco $e=(\Lambda,\Lambda')\in E T$ associamo il gruppo \[G_e=G_\Lambda\cap G_{\Lambda'},\] e l'inclusione $\psi_e:G_e\to G_{\Lambda'}$. Questi dati definiscono un albero di gruppi $\mathcal G(T)=(G_v,G_e,\psi_e)_{v\in \cli,e\in ET}$. La definizione dipende chiaramente dalla scelta dell'albero $T$. Tuttavia, vale:
\begin{thm}
	Il grafo di gruppi $\mathcal G(T)$ ha gruppo fondamentale relativo all'albero $T$ $$\pi_1(\mathcal G(T),T)=G_\Gamma.$$
	\begin{proof}
		Per una clique $\Lambda$, denotiamo con $x^\Lambda$ il generatore del gruppo $G_\Lambda$ corrispondente al vertice $x\in V\Lambda$. Poiché vale la proprietà di intersezione per $T$, si ha: \[\pi_1(\mathcal G(T),T)=\frac{\coprod_{\Lambda\in \cli} G_\Lambda}{x^\Lambda=x^{\Lambda'}:\ x\in V\Lambda\cap V\Lambda'}.\] Dunque il gruppo fondamentale ha presentazione \[\langle x^\Lambda:\ x\in\Lambda,\ \Lambda\in \cli\ \vert\ R(e)=1,\ x^\Lambda=x^{\Lambda'}:\ e\in E\Lambda,\ x\in V\Lambda\cap V\Lambda'\rangle, \] che è uguale a \[\langle x: x\in V\ \vert \ R(e):\ e\in E\rangle=G_\Gamma.\]
	\end{proof}
\end{thm}
Dato un grafo speciale $\Gamma=(V,E)$, sia $W_\Gamma$ il \textit{gruppo right-angled di Coxeter} (RACG) associato a $\ddot \Gamma$, i.e. \[W_\Gamma=\langle V\ \vert\ v^2,\ [v_1,v_2]:\ v\in V,\ \{v_1,v_2\}\in\ddot E\rangle. \]
Risulta quindi definito un epimorfismo naturale $$\pi=\pi_\Gamma:G_\Gamma\to W_\Gamma,$$ definito da $\pi(v)=v$ per ogni generatore $v\in V$. Il nucleo di $\pi$ contiene un sottogruppo isomorfo al RAAG classico associato a $\Gamma$. \\

Tramite le considerazioni fatte in questa sezione, si possono definire alberi di gruppi $\mathcal A$ e $\mathcal W$ con gruppi fondamentali rispettivamente $A_\Gamma$ e $W_\Gamma$. Denotiamo poi con $1$ l'albero di gruppi basato su $T$ con gruppi locali banali.
\begin{prop}
Si ha una successione esatta corta di grafi di gruppi \[1\to \mathcal A\to\mathcal G\to\mathcal W\to 1,\] cioè, $\forall e\in ET$, $v=t(e)$, il diagramma
\[
\begin{tikzcd}
1 \arrow[r] & A_e \arrow[r] \arrow[d] & G_e \arrow[r] \arrow[d] & W_e \arrow[r] \arrow[d] & 1 \\
1 \arrow[r] & A_v \arrow[r]           & G_v \arrow[r]           & W_v \arrow[r]           & 1
\end{tikzcd}\] commuta e le righe sono esatte.
\end{prop}

\bibliography{sample}{}

\begin{thebibliography}{Dav17}

\bibitem[AH99]{HB}
Martin R.~Bridson André~Haefliger.
\newblock {\em Metric spaces of non-positive curvature}.
\newblock Springer-Verlag Berlin Heidelberg, 1999.

\bibitem[Cas19]{C}
Alberto Cassella.
\newblock Hypercubical groups.
\newblock Master's thesis, Università degli Studi di Milano-Bicocca,
  2018-2019.

\bibitem[Dav17]{D}
Donald~M. Davis.
\newblock {\em n-dimensional Klein bottles}.
\newblock arXiv:1706.03704, 2017.

\bibitem[Dro87]{Dr}
Carl Droms.
\newblock Isomorphisms of graph groups.
\newblock {\em American Mathematical Society}, 100(3), July 1987.

\bibitem[Euc]{Eu}
Euclide.
\newblock Elementi.

\bibitem[Hat02]{hat}
Allen Hatcher.
\newblock {\em Algebraic Topology}.
\newblock Cambridge University Press, 2002.

\bibitem[JB93]{BP}
B.~Peyton J.R.S.~Blair.
\newblock {\em An Introduction to Chordal Graphs and Clique Trees}.
\newblock Springer, New York, NY, 1993.

\bibitem[Lan99]{L}
Serge Lang.
\newblock {\em Math Talks for Undergraduates}, chapter Bruhat-Tits Spaces.
\newblock Springer-Verlag New York, 1999.

\bibitem[Ser80]{S}
Jean-Pierre Serre.
\newblock {\em Trees, Translated from the French by John Stillwell}.
\newblock Springer-Verlag, 1980.

\end{thebibliography}
\bibliographystyle{alpha}

\end{document}